\documentclass[twocolumn]{autart}

\usepackage{graphicx}
\usepackage{amssymb}
\usepackage{amsmath}
\usepackage{bm}
\usepackage{tikz}

\usetikzlibrary{external}
\usepackage[round]{natbib}

\usepackage[hidelinks]{hyperref}
\hypersetup{
     colorlinks   = true,
     citecolor    = blue 
}
\begin{document}
    \begin{frontmatter}
        \title{Minimum-time interception of a moving target by a material point in a viscous medium\thanksref{footnoteinfo}} 
        
        \thanks[footnoteinfo]{This paper was not presented at any IFAC 
        meeting. Corresponding author M.~E.~Buzikov. Tel. +7 (966) 095-03-04.}
        
        \author[ICS]{Maksim E. Buzikov}\ead{me.buzikov@physics.msu.ru},
        \author[ICS]{Alina M. Mayer}\ead{mayer@ipu.ru}
        
        \address[ICS]{V.~A.~Trapeznikov Institute of Control Sciences of Russian Academy of Sciences, Moscow, Russia} 
                  
        
        \begin{abstract}
            In this paper we investigated a model that describes the motion of a material point in a viscous medium under a force that is arbitrary in direction, but limited in magnitude. This model was named "the isotropic rocket" in the early work of Rufus Isaacs. We obtained a parametric description of a reachable set for the isotropic rocket and solved a group of reachability problems for a final configuration that varies in a known time-dependent manner (moving target). To describe the reachable set, we obtained an explicit parametric form of its boundary and all its projections onto all subspaces of the state space. Convergent algorithms have been proposed for computing minimum-time interception in position and velocity spaces. Finally, we numerically investigated particular cases of minimum-time interception, validating the development of this study.
        \end{abstract}
    \end{frontmatter}

    \section{Introduction}
    
    The path-following problem occurs in complex missions for autonomous ground, aerial, and marine vehicles. This problem requires a given reference path as initial data (see \cite{Kai2019-pi,Yao2020-ye,Zheng2020-hx}). The reference path can be specified by the mission conditions, for example, for a terrain surveillance mission. In other scenarios, the mission only declares certain requirements of effectiveness or safety. The corresponding reference path should be computed onboard, and the computation should be time-rapid and memory-efficient. The main motivation of our study is to create an efficient and predictable method for computing the reference path for moving target interception missions. We will use the isotropic rocket as a simple model that produces twice differentiable reference paths (see \cite{Isaacs1965-nd}).

    The isotropic rocket is a linear model that describes the motion of a material point in a viscous medium. This model apparently received its name from the early work of \cite{Isaacs1955-bn} on the theory of differential games. Note that the name “isotropic rocket” was chosen rather dubiously, because the corresponding model is more suitable for constructing smooth trajectories than for describing the dynamics of a real missile (see~\cite{Manchester2002-sl}). Isaacs considered a pursuit-evasion game in which an isotropic rocket played the role of the pursuer and a plant with simple motion dynamics played the role of the evader. \cite{Isaacs1955-bn,Isaacs1965-nd} examined the problem of constructing feedback controls for the case of a nonzero drag coefficient, and the problem of constructing a barrier surface for a zero drag coefficient, including analytically verifiable conditions for avoiding capture. \cite{Bernhard1970-po,Bernhard1972-gd} used more rigorous constructions for the basic notions of differential game theory and clarified Isaacs' solution to the isotropic rocket game. Numerical investigations of the Bellman-Isaacs equation for an isotropic rocket game were performed by~\cite{Botkin2011-vp,Botkin2013-ox,Kumkov2014-bn}.

    The isotropic rocket appears in many variations of Isaacs' problem. \cite{Wong1967-lu} analyzed the game problem for the case of an isotropic rocket in the role of an evader, and the players were additionally affected by gravitational forces, but the drag coefficient was taken to be zero. For a wide class of linear differential games, \cite{Pontryagin1969-dt} formulated a sufficient condition for escape that becomes necessary in the corresponding surveillance evasion game. \cite{Lewin1989-ly} also analyzed the surveillance evasion game. \cite{Friedman1971-jn,Gutman1987-wz,Melikian1973-ee} considered an isotropic rocket game with a fixed termination time. Note that such problems have simpler solutions and explicit expressions for the value function. \cite{Abeysiriwardena2018-wv,Abeysiriwardena2019-yg} considered energy-efficient strategies in the game. Another studied variation of the game problem formulation is the case in which information regarding the trajectory of the evader movement is partially known. \cite{Selvakumar2015-pb,Selvakumar2018-np}  examined game problem formulation, implying the availability of only one piece of information regarding the evader trajectory.

    The isotropic rocket has been intensively studied for various optimal-control problems. \cite{Pontryagin1962-lt} obtained a minimum-time feedback control that transfers the isotropic rocket at a given point in the position space for a zero drag coefficient. A series of papers by L.D.~Akulenko is devoted to detailed studies of the minimum-time problems: reach a sphere in the position space (\cite{Akulenko1996-oq}), reach a given point in the position space with the initial value of velocity (\cite{Akulenko2003-xd}), return to the starting point with the required velocity (\cite{Akulenko2005-iz}), reach a point with a given velocity (\cite{Akulenko2007-aw}), reach an ellipsoid in the position space (\cite{Akulenko2018-ck}). The further studies by Akulenko et al., unlike the previous ones, are devoted to the case of a non-zero drag coefficient. \cite{Akulenko2007-rq} solved the minimum-time problem of transferring the isotropic rocket to a sphere in coordinate space (\cite{Vnuchkov1998-jc} analyzed the same problem). \cite{Akulenko2011-kr} considered the problem of the fastest transfer to a given point with a desired velocity value. In this study, the optimization problem is reduced to finding the roots of a system of real equations. The solution to these equations determines the open-loop control, calculates the corresponding trajectory of the isotropic rocket, and establishes the minimum transfer time. In our study, we supplemented analytical results of \cite{Akulenko2011-kr} by clarifying the singular cases of integration of the equations of motion for extremal control inputs. Note that the solution to the problem of transferring the isotropic rocket to a given state can be useful for interpolating a smooth curve passing through given points (see \cite{Tankasala2022-hn}).

    The current study is devoted to the problem of the minimum-time interception of a moving target by an isotropic rocket. It is assumed that the trajectory of the target is known a priori. In practice, target's trajectory is given a priori as a prediction, for example, when it is necessary to reveal the optimal order of visiting several targets (see \cite{Stieber2022-et}). Thus, an isotropic rocket must be transferred to a given set of states that can change over time in a known manner. This formulation includes cases of motion of an isotropic rocket under wind conditions considered by \cite{Akulenko2008-so,Bakolas2014-ze}. \cite{Bakolas2014-ze} obtained a real equation whose minimum non-negative root is the minimum time of interception. Intercepting a moving target here is equivalent to reaching a given position in a known changing wind. Using the minimum interception time, the optimal open-loop control can be explicitly computed. A similar problem was solved for wind, which depends linearly on the spatial location of an isotropic rocket (see \cite{Bakolas2014-jc,Bakolas2016-ry}). The result of these works is an algorithm for computing the value function for a defined time-space partitioning grid. This algorithm uses the ability to construct the boundary of a reachable set for a given model. The main contribution of our study is an always-convergent algorithm that computes a minimal time in the problem of minimum-time interception of a moving target by an isotropic rocket. In comparison with \cite{Bakolas2014-ze} results the proposed algorithm has proven convergence regardless of the initial approximation. We analytically describe the function of the distance from a given point to the projection of the reachable set to be able to use a framework of guaranteed computing (see \cite{Buzikov2022-rb}). This framework proposes an algorithm that computes the minimum interception time regardless of the initial guess, that is, it produces a minimal root of the \cite{Bakolas2014-ze} real equation. The design of the algorithm requires the ability to calculate the distance to the projection of the reachable set. We use the maximum principle to describe a boundary of the reachable set. The specific form of this boundary for the projections on position and velocity spaces permits to compute the distance function analytically. In addition, we provided a similar algorithm for the problem of reaching a desired time-varying velocity.

    \section{Problem formulation}
    
    In this section, we formulate the problem of interception of a moving target as an optimal control problem. Throughout this paper, we consider only Euclidean spaces with inner product and norm given by
    \begin{equation*}
        (\mathbf{a}, \mathbf{b}) = \sum_{i = 1}^k a_ib_i, \quad \lVert\mathbf{a}\rVert = \sqrt{(\mathbf{a},\mathbf{a})},
    \end{equation*}
    where $\mathbf{a} = \begin{bmatrix}a_1 & ... & a_k\end{bmatrix}^\top$, $\mathbf{b} = \begin{bmatrix}b_1 & ... & b_k\end{bmatrix}^\top$.

    \subsection{General problem}
    
    The state space of an isotropic rocket is $\mathcal{S} = \mathbb{R}^{2n}$, where $n$ is the dimension of the position space. The state of the isotropic rocket at time $t \in \mathbb{R}^+_0$ can be described using the position and velocity vectors:
    \begin{align*}
        &\boldsymbol{r}(t; \boldsymbol{u}) = 
        \begin{bmatrix}
            r^1(t; \boldsymbol{u}) & ... & r^n(t; \boldsymbol{u})
        \end{bmatrix}^\top,\\
        &\boldsymbol{v}(t; \boldsymbol{u}) = \begin{bmatrix}
            v^1(t; \boldsymbol{u}) & ... & v^n(t; \boldsymbol{u})
        \end{bmatrix}^\top.
    \end{align*}
    Here, $\boldsymbol{u} = \begin{bmatrix}u^1 & ... & u^n\end{bmatrix}^\top$ is a control input. Note that only the cases $n \in \{1, 2, 3\}$ are useful for practical purposes, but we will use arbitrary $n \in \mathbb{N}$ for the sake of completeness and uniformity of notation. We denote the state vector by the following way:
    \begin{equation*}
        \boldsymbol{s}(t; \boldsymbol{u}) = \begin{bmatrix}r^1(t; \boldsymbol{u}) & ... & r^n(t; \boldsymbol{u}) & v^1(t; \boldsymbol{u}) & ... & v^n(t; \boldsymbol{u})\end{bmatrix}^\top.
    \end{equation*}
    The isotropic rocket motion is described\footnote{For a model with arbitrary viscosity and maximal force magnitude, there exists a transformation of the time and space that makes these values units.} by
    \begin{equation}\label{eq:dyn_eq}
        \left\{
        \begin{aligned}
            &\dot{\boldsymbol{r}} = \boldsymbol{v};\\
            &\dot{\boldsymbol{v}} = \boldsymbol{u} - \boldsymbol{v}.
        \end{aligned}
        \right.
    \end{equation}
    The set of admissible control inputs is denoted by $\mathcal{A}$. These are measurable functions and $\boldsymbol{u}(t) \in \mathcal{U} = \{\mathbf{u} \in \mathbb{R}^n:\: \lVert\mathbf{u}\rVert \leq 1\}$ for any $\boldsymbol{u} \in \mathcal{A}$. The coordinate axes are aligned and rotated such that the initial position and velocity are the following:
    \begin{equation*}
        \boldsymbol{r}(0; \boldsymbol{u}) = \boldsymbol{0}, \quad \boldsymbol{v}(0; \boldsymbol{u}) = \begin{bmatrix}v_0 & 0 & ... & 0\end{bmatrix}^\top.
    \end{equation*}
    Here, $v_0 \in [0, 1)$ is an initial speed. Throughout the paper, we suppose that $\boldsymbol{r}(\cdot; \boldsymbol{u})$, $\boldsymbol{v}(\cdot; \boldsymbol{u})$ are absolutely continuous solutions of system~\eqref{eq:dyn_eq} for a fixed $\boldsymbol{u} \in \mathcal{A}$.

    Let $\mathcal{H}$ be a subspace of $\mathcal{S}$. $\mathcal{H}$ denotes the space of the state-vector components that are important for interception. The trajectory of the target is given by the function $\boldsymbol{h}_T \in \mathrm{Lip}_v(\mathbb{R}^+_0, \mathcal{H})$\footnote{$\mathrm{Lip}_v(\mathcal{F}, \mathcal{G})$ denotes a set of $v$-Lipschitz continuous functions. If $\boldsymbol{y} \in \mathrm{Lip}_v(\mathcal{F}, \mathcal{G})$, then $\boldsymbol{y}: \mathcal{F} \to \mathcal{G}$ and $\lVert \boldsymbol{y}(\mathbf{x}_2) - \boldsymbol{y}(\mathbf{x}_1) \rVert_\mathcal{G} \leq v \lVert \mathbf{x}_2 - \mathbf{x}_1 \rVert_\mathcal{F}$ for all $\mathbf{x}_1, \mathbf{x}_2 \in \mathcal{F}$.}, where $v \in \mathbb{R}^+_0$. The state vector components that are important for interception are denoted as $\boldsymbol{h}(t; \boldsymbol{u})$. We define the general problem of the minimum time interception of a moving target by finding a control input $\boldsymbol{u} \in \mathcal{A}$ that sets the minimum value of the following functional:
    \begin{equation*}
        J[\boldsymbol{u}; \boldsymbol{h}_T] \overset{\mathrm{def}}{=} \min \left\{t \in \mathbb{R}^+_0:\: \lVert\boldsymbol{h}(t; \boldsymbol{u}) - \boldsymbol{h}_T(t)\rVert \leq \ell\right\}.
    \end{equation*}
    Here, $\ell \in \mathbb{R}^+_0$ denotes a capture radius. The minimum-time required for an interception is defined by
    \begin{equation*}
        T^*[\boldsymbol{h}_T] \overset{\mathrm{def}}{=} \inf_{\boldsymbol{u} \in \mathcal{A}} J[\boldsymbol{u}; \boldsymbol{h}_T].
    \end{equation*} 

    \subsection{Two specific problems}

    We present two specific problems, which are special cases of the general problem of minimum time interception. The first problem involves intercepting a moving point in the position space.

    \begin{prob}\label{prob:1}
        Find the optimal control input if $\mathcal{H} = \mathbb{R}^n$, $\boldsymbol{h}(\cdot; \boldsymbol{u}) = \boldsymbol{r}(\cdot; \boldsymbol{u})$.
    \end{prob}

\begin{figure}
        \begin{center}
            \includegraphics[width=0.35\textwidth]{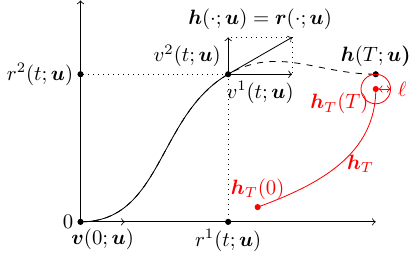}
            \caption{Transferring to a moving point in the position space for $n = 2$. The red line is a path of the moving point. The black line is a path of the isotropic rocket. The isotropic rocket captured the moving point at $T$}
            \label{fig:positional_interception}
        \end{center}
    \end{figure}

    The two-dimensional case of the position space for the problem is illustrated in Fig.~\ref{fig:positional_interception}. The trajectory of the target describes the motion of point $\boldsymbol{h}_T(t)$. The isotropic rocket must approach the point by distance $\ell$ as quickly as possible.
    
    The second problem is to reach a time-varying velocity.

    \begin{prob}\label{prob:2}
        Find the optimal control input if $\mathcal{H} = \mathbb{R}^n$, $\boldsymbol{h}(\cdot; \boldsymbol{u}) = \boldsymbol{v}(\cdot; \boldsymbol{u})$.
    \end{prob}

\begin{figure}
        \begin{center}
            \includegraphics[width=0.4\textwidth]{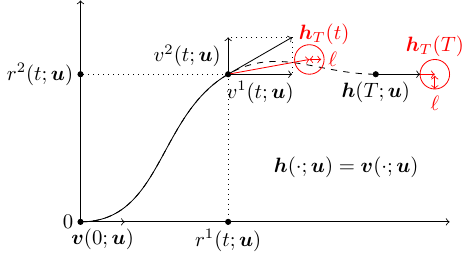}
            \caption{Reaching a desired time-varying velocity for $n = 2$. The red vectors are the desired velocity values at $t$ and $T$. The isotropic rocket reached the desired value of the velocity vector with an error $\ell$ at $T$}
            \label{fig:reach_velocity}
        \end{center}
    \end{figure}

    The two-dimensional case for the problem is illustrated in Fig.~\ref{fig:reach_velocity}. The trajectory of the target describes the time-varying velocity $\boldsymbol{h}_T(t)$. The isotropic rocket must reach the desired velocity with an error $\ell$ as quickly as possible.
    
    \section{The reachable set}
    
    In this section, we describe a reachable set of isotropic rocket. The reachable set consists of points in the state space that can be reached using the admissible control inputs: $\mathcal{R}(t) \overset{\mathrm{def}}{=} \left\{\boldsymbol{s}(t; \boldsymbol{u}):\: \boldsymbol{u} \in \mathcal{A}\right\}$. The reachable set is compact and convex for any $t \in \mathbb{R}^+_0$ and the set-valued mapping $\mathcal{R}$ is continuous (see \cite{Lee1967-vf}).

    \subsection{Extremal control inputs}

    According to \cite{Lee1967-vf}, maximal control inputs are extremal in linear processes theory. Thus, maximal control inputs are enough to reach any point of the boundary of the reachable set $\partial\mathcal{R}(t)$. Let
    \begin{align*}
        &\boldsymbol{\lambda}(t; T, \mathbf{p}_T) = \begin{bmatrix}\lambda_1(t; T, \mathbf{p}_T) & ... &\lambda_n(t; T, \mathbf{p}_T)\end{bmatrix},\\
        &\boldsymbol{\eta}(t; T, \mathbf{p}_T) = \begin{bmatrix}\eta_1(t; T, \mathbf{p}_T) & ... &\eta_n(t; T, \mathbf{p}_T)\end{bmatrix}
    \end{align*}
    be adjoint vectors for $\boldsymbol{r}(t; \boldsymbol{u})$, $\boldsymbol{v}(t; \boldsymbol{u})$. Here,
    \begin{align*}
        &\mathbf{p}_T = \begin{bmatrix}\lambda_1^T & ... & \lambda_n^T & \eta_1^T & ... & \eta_n^T\end{bmatrix},\\
        &\boldsymbol{\lambda}_T = \begin{bmatrix}\lambda_1^T & ... & \lambda_n^T\end{bmatrix}, \quad \boldsymbol{\eta}_T = \begin{bmatrix}\eta_1^T & ... & \eta_n^T\end{bmatrix}
    \end{align*}
    are the terminal values of adjoint vectors ($\boldsymbol{\lambda}(T; T, \mathbf{p}_T) = \boldsymbol{\lambda}_T$, $\boldsymbol{\eta}(T; T, \mathbf{p}_T) = \boldsymbol{\eta}_T$) at some terminal time moment $T \in \mathbb{R}^+_0$. The adjoint differential system $\dot{\boldsymbol{\lambda}} = \boldsymbol{0}$, $\dot{\boldsymbol{\eta}} = \boldsymbol{\eta} - \boldsymbol{\lambda}$ has an explicit form of solution:
    \begin{equation}\label{eq:lambda_eta}
        \boldsymbol{\lambda}(t; T, \mathbf{p}_T) = \boldsymbol{\lambda}_T, \quad \boldsymbol{\eta}(t; T, \mathbf{p}_T) = \boldsymbol{\lambda}_T + (\boldsymbol{\eta}_T - \boldsymbol{\lambda}_T)e^{t - T}.
    \end{equation}
    Let $\boldsymbol{u}_E(\cdot; T, \mathbf{p}_T) \in \mathcal{A}$ be an extremal control input satisfying the maximum principle:
    \begin{equation}\label{eq:max_princ}
        \boldsymbol{\eta}(t; T, \mathbf{p}_T)\boldsymbol{u}_E(t; T, \mathbf{p}_T) \overset{\mathrm{a.e.}}{=} \max_{\mathbf{u} \in \mathcal{U}} \boldsymbol{\eta}(t; T, \mathbf{p}_T)\mathbf{u}.
    \end{equation}
    Before obtaining an explicit form for $\boldsymbol{u}_E(\cdot; T, \mathbf{p}_T)$, we obtain conditions for $\boldsymbol{\eta}(t; T, \mathbf{p}_T) = \boldsymbol{0}$.
    \begin{lem}\label{lem:eta_zero}
        There exists $t \in (0, T]$ such that $\boldsymbol{\eta}(t; T, \mathbf{p}_T) = \boldsymbol{0}$ if and only if
        \begin{equation}\label{eq:eta_zero_cond}
            \lVert\boldsymbol{\lambda}_T\rVert\lVert\boldsymbol{\eta}_T\rVert = -(\boldsymbol{\lambda}_T, \boldsymbol{\eta}_T), \quad e^{-T} < \frac{\lVert\boldsymbol{\lambda}_T\rVert}{\lVert\boldsymbol{\lambda}_T\rVert + \lVert\boldsymbol{\eta}_T\rVert}.
        \end{equation}
        This value $t \in (0, T]$ is unique and
        \begin{equation}\label{eq:switch_time}
            t = \theta(T, \mathbf{p}_T) \overset{\mathrm{def}}{=} T + \ln{\frac{\lVert\boldsymbol{\lambda}_T\rVert}{\lVert\boldsymbol{\lambda}_T\rVert + \lVert\boldsymbol{\eta}_T\rVert}}.
        \end{equation}
    \end{lem}
    \begin{pf}
        At first, we prove that if there exists $t \in (0, T]$, such that
        \begin{equation*}
            \boldsymbol{\eta}(t; T, \mathbf{p}_T) = \boldsymbol{\lambda}_T + (\boldsymbol{\eta}_T - \boldsymbol{\lambda}_T)e^{t - T} = \boldsymbol{0},
        \end{equation*}
        then conditions~\eqref{eq:eta_zero_cond} are satisfied. It is obvious, that $\boldsymbol{\lambda}_T$, $\boldsymbol{\eta}_T$ are collinear and $\boldsymbol{\eta}_T = \boldsymbol{\lambda}_T(1 - e^{T - t})$. Using $1 - e^{T - t} \leq 0$, we obtain
        \begin{equation*}
            (\boldsymbol{\lambda}_T, \boldsymbol{\eta}_T) =
            -\lVert\boldsymbol{\lambda}_T\rVert\lVert\boldsymbol{\lambda}_T(e^{T - t} - 1)\rVert = -\lVert\boldsymbol{\lambda}_T\rVert\lVert\boldsymbol{\eta}_T\rVert.
        \end{equation*}
        Combining $\lVert\boldsymbol{\eta}_T\rVert = \lVert\boldsymbol{\lambda}_T\rVert(e^{T - t} - 1)$ and $e^{-t} < 1$, we yield
        \begin{equation*}
            e^{-T}(\lVert\boldsymbol{\eta}_T\rVert + \lVert\boldsymbol{\lambda}_T\rVert) = \lVert\boldsymbol{\lambda}_T\rVert e^{-t} < \lVert\boldsymbol{\lambda}_T\rVert.
        \end{equation*}
        That is, conditions~\eqref{eq:eta_zero_cond} are satisfied. Now we demonstrate their sufficiency. We prove that the value $t$ given by~\eqref{eq:switch_time} satisfying~\eqref{eq:eta_zero_cond} is the root of $\boldsymbol{\eta}(t; T, \mathbf{p}_T) = \boldsymbol{0}$ and $0 < t \leq T$. The inequality $t \leq T$ is obvious, and $0 < t$ follows from
        \begin{equation*}
            e^{-T} < \frac{\lVert\boldsymbol{\lambda}_T\rVert}{\lVert\boldsymbol{\lambda}_T\rVert + \lVert\boldsymbol{\eta}_T\rVert} = e^{t - T}.
        \end{equation*}
        $\boldsymbol{\lambda}_T \neq \boldsymbol{0}$, because if $\boldsymbol{\lambda}_T = \boldsymbol{0}$, then the second inequality from~\eqref{eq:eta_zero_cond} is violated. The first condition from~\eqref{eq:eta_zero_cond} says that the vectors $\boldsymbol{\lambda}_T$, $\boldsymbol{\eta}_T$ are collinear and directed in opposite directions, i.e., $\boldsymbol{\eta}_T = -\boldsymbol{\lambda}_T\lVert\boldsymbol{\eta}_T\rVert/\lVert\boldsymbol{\lambda}_T\rVert$. The equation~\eqref{eq:switch_time} gives $\lVert\boldsymbol{\eta}_T\rVert/\lVert\boldsymbol{\lambda}_T\rVert = e^{T - t} - 1$. It implies
        \begin{equation*}
            \boldsymbol{\eta}_T = -\boldsymbol{\lambda}_T\frac{\lVert\boldsymbol{\eta}_T\rVert}{\lVert\boldsymbol{\lambda}_T\rVert} = \boldsymbol{\lambda}_T(1 - e^{T - t}),
        \end{equation*}
        which proves $\boldsymbol{\eta}(t; T, \mathbf{p}_T) = \boldsymbol{0}$. Now we prove the uniqueness of $t$ defined by~\eqref{eq:switch_time}. Consider
        \begin{multline*}
            \lVert\boldsymbol{\eta}(t; T, \mathbf{p}_T)\rVert^2 = \lVert\boldsymbol{\eta}_T - \boldsymbol{\lambda}_T\rVert^2\left(e^{t - T}\right)^2\\
            + 2(\boldsymbol{\lambda}_T, \boldsymbol{\eta}_T - \boldsymbol{\lambda}_T)e^{t - T} + \lVert\boldsymbol{\lambda}_T\rVert^2 = 0.
        \end{multline*}
        The equation is quadratic with respect to $e^{t - T}$, and its discriminant is calculated as follows:
        \begin{equation*}
            D = 4(\boldsymbol{\lambda}_T, \boldsymbol{\eta}_T - \boldsymbol{\lambda}_T)^2 - 4\lVert\boldsymbol{\eta}_T - \boldsymbol{\lambda}_T\rVert^2\lVert\boldsymbol{\lambda}_T\rVert^2 = 0.
        \end{equation*}
        The last equality is true due to the first condition from~\eqref{eq:eta_zero_cond}. The monotonicity of $e^{t - T}$ over $t$ and the uniqueness of the solution of the quadratic equation for $D = 0$ prove the uniqueness of the root~\eqref{eq:switch_time}.\qed
    \end{pf}

    \begin{lem}\label{lem:ext_cont}
        Any point of the boundary of reachable set $\partial\mathcal{R}(T)$ can be visited using the control inputs given by
        \begin{equation}\label{eq:uE}
            \boldsymbol{u}_E(t; T, \mathbf{p}_T) \overset{\mathrm{a.e.}}{=} \frac{\boldsymbol{\eta}(t; T, \mathbf{p}_T)}{\lVert\boldsymbol{\eta}(t; T, \mathbf{p}_T)\rVert},
        \end{equation}
        where $\mathbf{p}_T \in \mathbb{R}^{2n}$, $\lVert \mathbf{p}_T \rVert = 1$. Moreover, these control inputs lead to $\partial\mathcal{R}(T)$ for any $\mathbf{p}_T \in \mathbb{R}^{2n}$.
    \end{lem}
    \begin{pf}
        If $t \in [0, T]$ and $\boldsymbol{\eta}(t; T, \mathbf{p}_T) \neq \boldsymbol{0}$, then~\eqref{eq:max_princ} yields~\eqref{eq:uE}. Due to Lemma~\ref{lem:eta_zero}, $\boldsymbol{\eta}(t; T, \mathbf{p}_T) \neq \boldsymbol{0}$ for almost all $t \in [0, T]$. It implies that the expression~\eqref{eq:uE} defines extreme control inputs almost everywhere. If $\lVert\mathbf{p}_T\rVert \neq 1$, then we can divide the numerator and denominator of the expression~\eqref{eq:uE} by $\lVert\mathbf{p}_T\rVert$ and 
        \begin{equation*}
            \boldsymbol{u}_E(t; T, \mathbf{p}_T) \overset{\text{a.e.}}{=} \frac{\boldsymbol{\eta}(t; T, \mathbf{p}_T)}{\lVert\boldsymbol{\eta}(t; T, \mathbf{p}_T)\rVert} = \frac{\boldsymbol{\eta}(t; T, \mathbf{p}_T/\lVert\mathbf{p}_T\rVert)}{\lVert\boldsymbol{\eta}(t; T, \mathbf{p}_T/\lVert\mathbf{p}_T\rVert)\rVert},
        \end{equation*}
        for $\lVert\mathbf{p}_T\rVert \neq 0$. That is, $\boldsymbol{u}_E(t; T, \mathbf{p}_T) = \boldsymbol{u}_E(t; T, \mathbf{p}_T/\lVert\mathbf{p}_T\rVert)$. Theorem 2 of~\cite[p.~73]{Lee1967-vf} guarantees that there are enough maximal control inputs to cover the boundary of the reachable set. Moreover, there are no maximal control inputs leading inside the reachable set.\qed
    \end{pf}

    \subsection{Extremal trajectories}

    We can now obtain a description of the extremal trajectories leading to $\partial\mathcal{R}(t)$. Integrating the equations of motion~\eqref{eq:dyn_eq} with extremal control inputs~\eqref{eq:uE} gives
    \begin{multline*}
        \boldsymbol{v}_E(t; T, \mathbf{p}_T) \overset{\mathrm{def}}{=} \boldsymbol{v}(t; \boldsymbol{u}_E(\cdot; T, \mathbf{p}_T))\\
        = \boldsymbol{v}_0 e^{-t} + \int\limits_0^t \boldsymbol{u}_E(s; T, \mathbf{p}_T)e^{s - t}\mathrm{d}s;
    \end{multline*}
    \begin{multline*}
        \boldsymbol{r}_E(t; T, \mathbf{p}_T) \overset{\mathrm{def}}{=} \boldsymbol{r}(t; \boldsymbol{u}_E(\cdot; T, \mathbf{p}_T)) \\
        = \boldsymbol{v}_0 - \boldsymbol{v}_E(t; T, \mathbf{p}_T) + \int\limits_0^t \boldsymbol{u}_E(s; T, \mathbf{p}_T)\mathrm{d}s.
    \end{multline*}
    Integrals from the right side of the expressions have singularities of the integrand function only if the conditions of Lemma~\ref{lem:eta_zero} are satisfied. Appendix~\ref{ap:ex_traj} analyzes each of four possible cases of behavior of the integrand function depending on the terminal value of the adjoint vector $\mathbf{p}_T$. The following expressions are valid for the trajectories leading to the boundary of the reachable set:
    \begin{multline}\label{eq:vE}
        \boldsymbol{v}_E(t; T, \mathbf{p}_T) \\ =
        \begin{cases}
            \text{\eqref{eq:vE_1}}, \: \lVert\boldsymbol{\lambda}_T\rVert\lVert\boldsymbol{\eta}_T\rVert > |(\boldsymbol{\lambda}_T, \boldsymbol{\eta}_T)|;\\
            \text{\eqref{eq:vE_2}}, \: \lVert\boldsymbol{\lambda}_T\rVert\lVert\boldsymbol{\eta}_T\rVert = -(\boldsymbol{\lambda}_T, \boldsymbol{\eta}_T), \: \boldsymbol{\lambda}_T \neq \boldsymbol{0};\\
            \text{\eqref{eq:vE_3}}, \: \lVert\boldsymbol{\lambda}_T\rVert\lVert\boldsymbol{\eta}_T\rVert = (\boldsymbol{\lambda}_T, \boldsymbol{\eta}_T), \: \boldsymbol{\lambda}_T \neq \boldsymbol{0};\\
            \text{\eqref{eq:vE_4}}, \: \boldsymbol{\lambda}_T = \boldsymbol{0}.
        \end{cases}
    \end{multline}
    \begin{multline}\label{eq:rE}
        \boldsymbol{r}_E(t; T, \mathbf{p}_T) \\ =
        \begin{cases}
            \text{\eqref{eq:rE_1}}, \: \lVert\boldsymbol{\lambda}_T\rVert\lVert\boldsymbol{\eta}_T\rVert > |(\boldsymbol{\lambda}_T, \boldsymbol{\eta}_T)|;\\
            \text{\eqref{eq:rE_2}}, \: \lVert\boldsymbol{\lambda}_T\rVert\lVert\boldsymbol{\eta}_T\rVert = -(\boldsymbol{\lambda}_T, \boldsymbol{\eta}_T), \: \boldsymbol{\lambda}_T \neq \boldsymbol{0};\\
            \text{\eqref{eq:rE_3}}, \: \lVert\boldsymbol{\lambda}_T\rVert\lVert\boldsymbol{\eta}_T\rVert = +(\boldsymbol{\lambda}_T, \boldsymbol{\eta}_T), \: \boldsymbol{\lambda}_T \neq \boldsymbol{0};\\
            \text{\eqref{eq:rE_4}}, \: \boldsymbol{\lambda}_T = \boldsymbol{0}.
        \end{cases}
    \end{multline}
    Note, that Eqs.~2.3--2.9 of~\cite{Akulenko2011-kr} are derived only for $\lVert\boldsymbol{\lambda}_T\rVert\lVert\boldsymbol{\eta}_T\rVert > |(\boldsymbol{\lambda}_T, \boldsymbol{\eta}_T)|$. Thus,~\eqref{eq:vE},~\eqref{eq:rE} clarify results of~\cite{Akulenko2011-kr}.
    
\begin{figure*}
        \centering
        \includegraphics[width=0.20\textwidth]{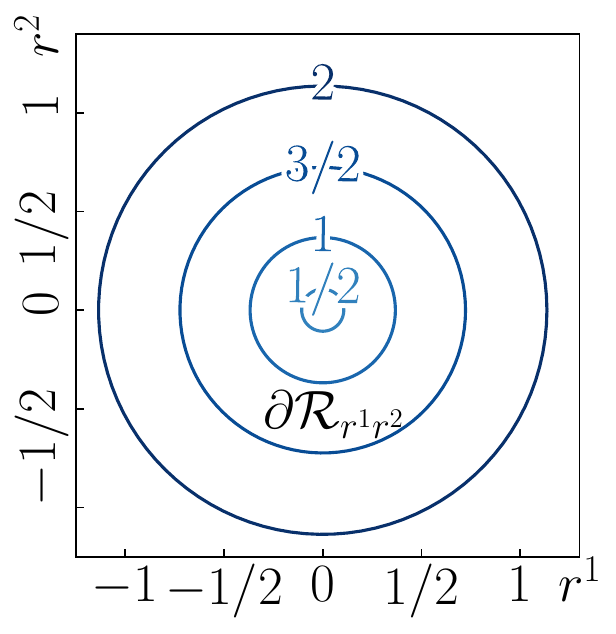}
        \includegraphics[width=0.255\textwidth]{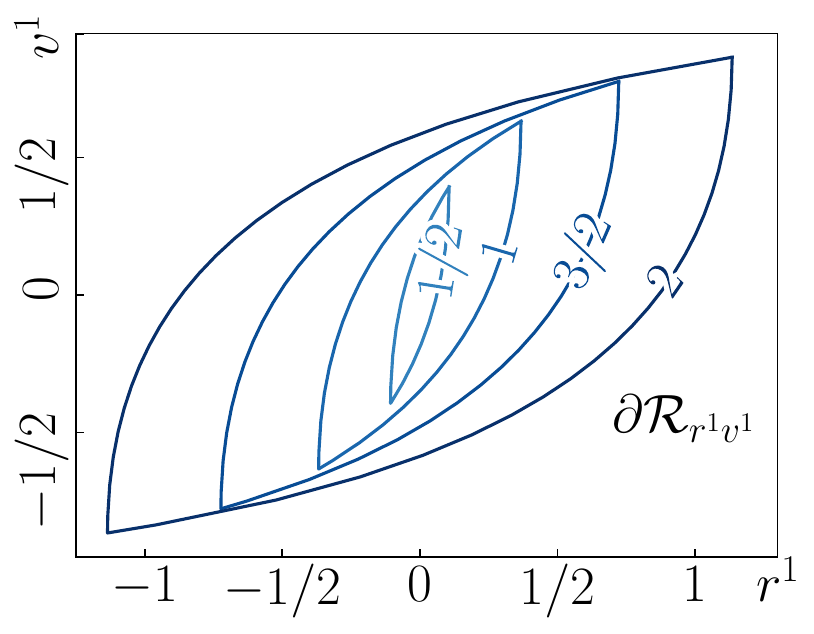}
        \includegraphics[width=0.20\textwidth]{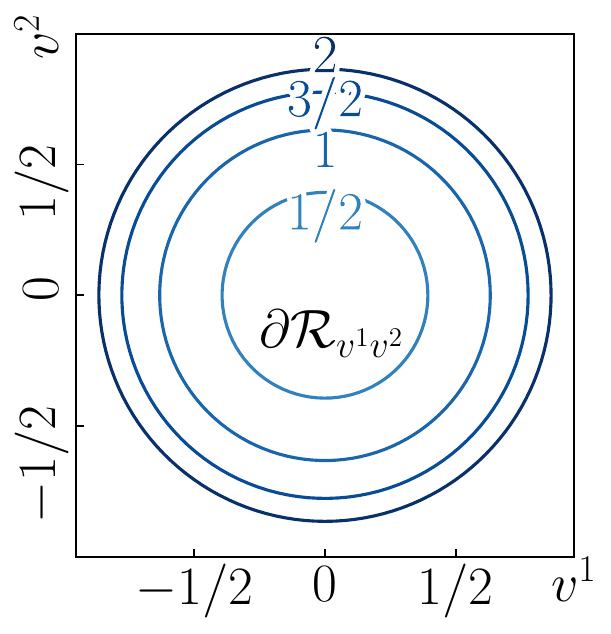}
        \includegraphics[width=0.255\textwidth]{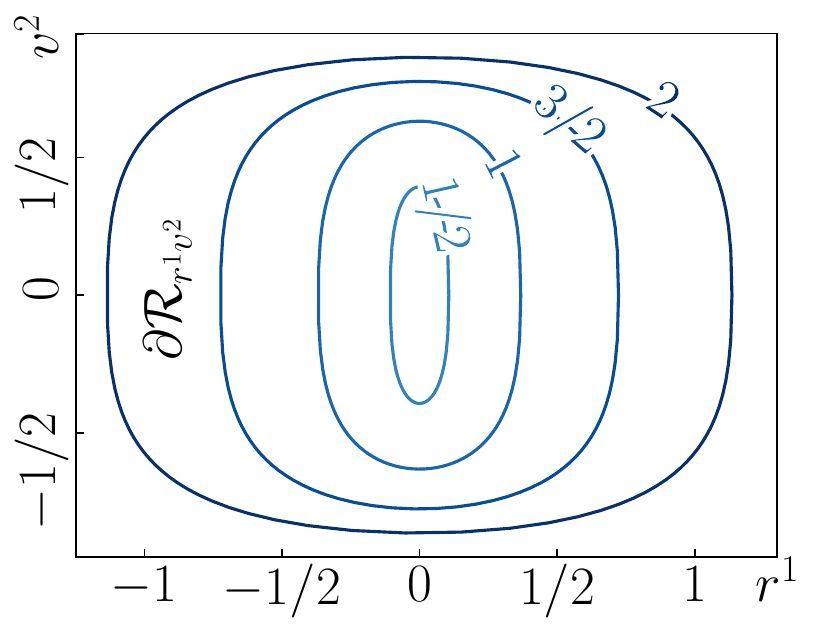}
        \includegraphics[width=0.20\textwidth]{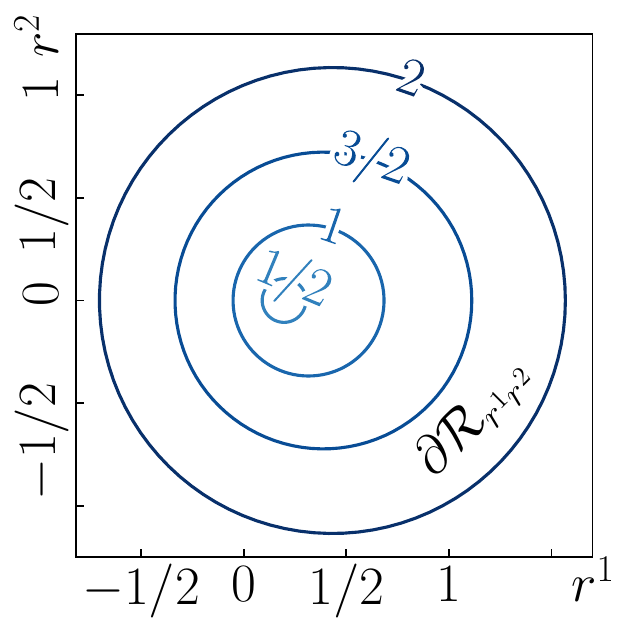}
        \includegraphics[width=0.255\textwidth]{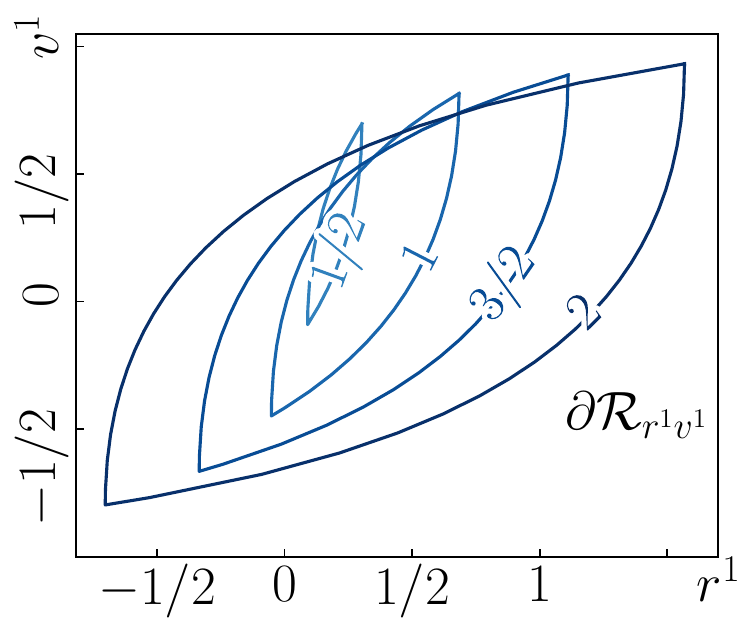}
        \includegraphics[width=0.20\textwidth]{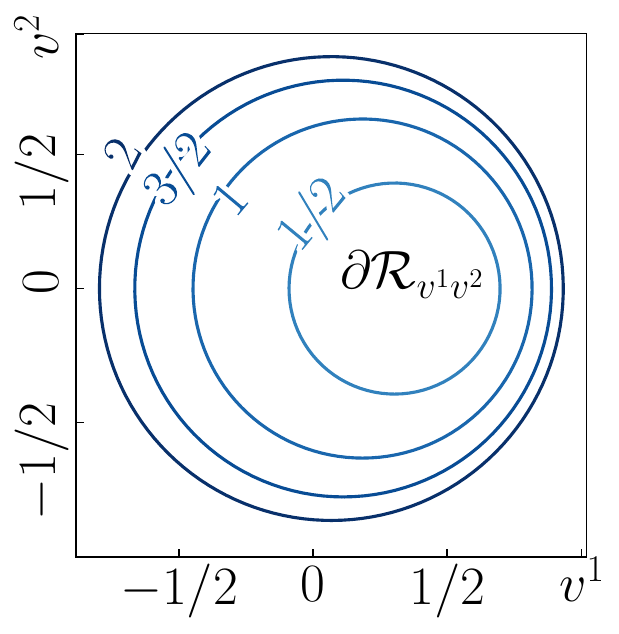}
        \includegraphics[width=0.255\textwidth]{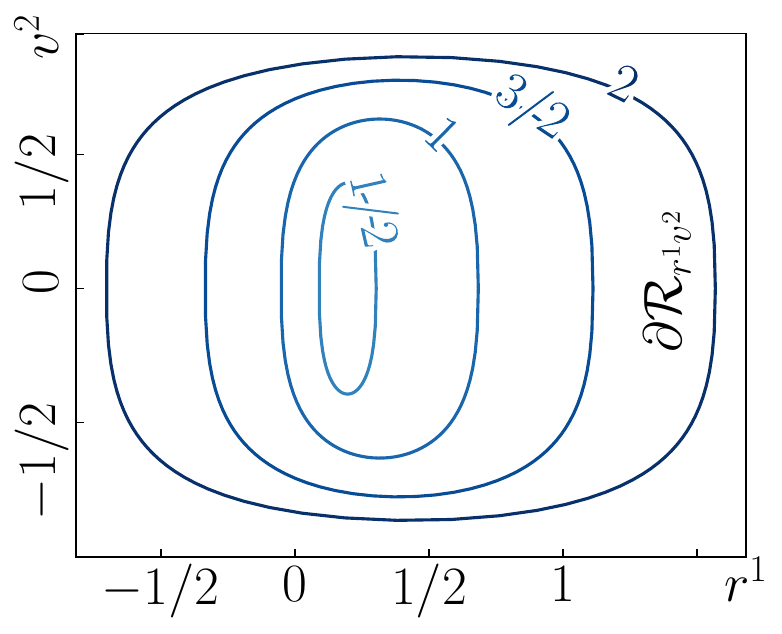}
        \caption{Boundaries of projections of the reachable set onto two-dimensional subspaces $r^1r^2$, $r^1v^1$, $v^1v^2$, $r^1v^2$ for different moments of time $t \in \{1/2, 1, 3/2, 2\}$ and initial velocities $v_0 = 0$ (upper row), $v_0 = 1/2$ (lower row). The state-space is four-dimensional ($n = 2$).}
        \label{fig:2d}
    \end{figure*}

    Define
    \begin{multline*}
        \boldsymbol{s}_E(t; T, \mathbf{p}_T) \overset{\mathrm{def}}{=} \boldsymbol{s}(t; \boldsymbol{u}_E(\cdot; T, \mathbf{p}_T)) \\
        = \mathrm{stack}(\boldsymbol{r}_E(t; T, \mathbf{p}_T), \boldsymbol{v}_E(t; T, \mathbf{p}_T)).
    \end{multline*}
    The boundary of the reachable set $\partial\mathcal{R}(T)$ consists of the endpoints of extremal trajectories. Hence, the explicit expression for this surface is given by
    \begin{equation}\label{eq:dR}
        \partial\mathcal{R}(T) = \{\boldsymbol{s}_E(T; T, \mathbf{p}_T):\: \mathbf{p}_T \in \mathbb{R}^{2n}, \: \lVert\mathbf{p}_T\rVert = 1\}.
    \end{equation}
    For the following purposes we will use special notations for the projections of the reachable set onto different subspaces of the state space. Let the numbers $i_1, ..., i_k \in \{1, ..., n\}$, $j_1, ..., j_m \in \{1, ..., n\}$, where $k, m \in \{0, ..., n\}$, be such that $i_1 < i_2 < .... < i_k$ and $j_1 < j_2 < ... < j_m$. We denote the projection of the reachable set at time $t$ onto a subspace of the state space with components $r^{i_1}...r^{i_k}v^{j_1}...v^{j_m}$ as $\mathcal{R}_{r^{i_1}...r^{i_k}v^{j_1}...v^{j_m}}(t)$ (see Fig.~\ref{fig:2d}).
    
    \section{Problem solution}

    In this section we propose always convergent algorithms for Problems~\ref{prob:1},~\ref{prob:2}. We use a description of the reachable set obtained in the previous section. For brevity, we will use the following notation:
    \begin{equation*}
        \mathcal{R}_{\boldsymbol{r}}(t) \overset{\mathrm{def}}{=} \mathcal{R}_{r^1...r^n}(t), \quad \mathcal{R}_{\boldsymbol{v}}(t) \overset{\mathrm{def}}{=} \mathcal{R}_{v^1...v^n}(t).
    \end{equation*}
    Here, $\mathcal{R}_{\boldsymbol{r}}(t)$ is a projection of the reachable set on the position space and $\mathcal{R}_{\boldsymbol{v}}(t)$ is a projection on the velocity space. We denote the distances to these sets as follows:
    \begin{align*}
        &\rho_{\mathcal{R}_{\boldsymbol{r}}}(t, \mathbf{h}) \overset{\mathrm{def}}{=} \min_{\tilde{\mathbf{h}} \in \mathcal{R}_{\boldsymbol{r}}(t)}\lVert\mathbf{h} - \tilde{\mathbf{h}}\rVert,\\
        &\rho_{\mathcal{R}_{\boldsymbol{v}}}(t, \mathbf{h}) \overset{\mathrm{def}}{=} \min_{\tilde{\mathbf{h}} \in \mathcal{R}_{\boldsymbol{v}}(t)}\lVert\mathbf{h} - \tilde{\mathbf{h}}\rVert.
    \end{align*}
     We use $\rho_{\mathcal{R}_{\boldsymbol{r}}}$ for solution of minimum-time interception problem in coordinate space (Problem~\ref{prob:1}) and $\rho_{\mathcal{R}_{\boldsymbol{v}}}$ is used for the problem of reaching a desired time-varying velocity (Problem~\ref{prob:2}).

    It follows from the transversality conditions that for Problem~\ref{prob:1} we have $\boldsymbol{\eta}_T = \boldsymbol{0}$, and for Problem~\ref{prob:2} we have $\boldsymbol{\lambda}_T = \boldsymbol{0}$. Therefore, we derive from~\eqref{eq:dR},~\eqref{eq:vE}, and~\eqref{eq:rE} that
    \begin{multline*}
        \partial\mathcal{R}_{\boldsymbol{r}}(t) = \{\boldsymbol{v}_0(1 - e^{-t}) + \frac{\boldsymbol{\lambda}_T}{\lVert\boldsymbol{\lambda}_T\rVert}(t - 1 + e^{-t}):\\
        \boldsymbol{\lambda}_T \in \mathbb{R}^n, \: \lVert\boldsymbol{\lambda}_T\rVert = 1\},
    \end{multline*}
   \begin{multline*}
        \partial\mathcal{R}_{\boldsymbol{v}}(t) = \{\boldsymbol{v}_0 e^{-t} + \frac{\boldsymbol{\eta}_T}{\lVert\boldsymbol{\eta}_T\rVert}(1 - e^{-t}):\\
        \boldsymbol{\eta}_T \in \mathbb{R}^n, \: \lVert\boldsymbol{\eta}_T\rVert = 1\}.
    \end{multline*} 
    Due to the convexity of the reachable set, we can conclude that the sets $\mathcal{R}_{\boldsymbol{r}}(t)$, $\mathcal{R}_{\boldsymbol{v}}(t)$ are balls with centers at points $\boldsymbol{v}_0(1 - e^{-t})$, $\boldsymbol{v}_0 e^{-t}$ and radii $t - 1 + e^{-t}$, $1 - e^{-t}$, respectively. Thus, the distances to these sets can be calculated as follows:
    \begin{align*}
        &\rho_{\mathcal{R}_{\boldsymbol{r}}}(t, \mathbf{h}) = \max(0, \lVert\mathbf{h} - \boldsymbol{v}_0(1 - e^{-t})\rVert - (t - 1 + e^{-t})),\\
        &\rho_{\mathcal{R}_{\boldsymbol{v}}}(t, \mathbf{h}) = \max(0, \lVert\mathbf{h} - \boldsymbol{v}_0 e^{-t}\rVert - (1 - e^{-t})).
    \end{align*}
    According to Lemma~3.6 of \cite{Buzikov2022-rb}, the sequence
    \begin{equation}\label{eq:seq_tau}
        t_i = \tau_\rho(t_{i - 1}, \boldsymbol{h}_T(t_{i - 1})), \quad t_0 = 0,
    \end{equation}
    converges to minimum interception time $T^*[\boldsymbol{h}_T]$. Here, $\tau_{\rho}$ is the simple universal lower estimator:
    \begin{equation*}
        \tau_\rho(t, \mathbf{h}) \overset{\mathrm{def}}{=}
        \begin{cases}
            t + \frac{\rho(t, \mathbf{h}) - \ell}{v_{\max} + v}, \quad \rho(t, \mathbf{h}) > \ell;\\
            t, \quad \rho(t, \mathbf{h}) \leq \ell.
        \end{cases}
    \end{equation*}
    Here, $\rho_{\mathcal{R}_{\boldsymbol{r}}}$ or $\rho_{\mathcal{R}_{\boldsymbol{v}}}$ is substituted in place of the function $\rho$, depending on the problem being solved. Lemma 3.6 also requires $\lVert\dot{\boldsymbol{h}}(t; \boldsymbol{u})\rVert \leq v_{\max}$, where $v_{\max} \in \mathbb{R}^+_0$. In case of Problem~\ref{prob:1}, we have
    \begin{equation*}
        \lVert\dot{\boldsymbol{h}}(t; \boldsymbol{u})\rVert = \lVert\dot{\boldsymbol{r}}(t; \boldsymbol{u})\rVert = \lVert\boldsymbol{v}(t; \boldsymbol{u})\rVert \leq 1 = v_{\max}.
    \end{equation*}
    For the Problem~\ref{prob:2}, we obtain
    \begin{multline*}
        \lVert\dot{\boldsymbol{h}}(t; \boldsymbol{u})\rVert = \lVert\dot{\boldsymbol{v}}(t; \boldsymbol{u})\rVert = \lVert\boldsymbol{u}(t) - \boldsymbol{v}(t; \boldsymbol{u})\rVert\\
        \leq \lVert\boldsymbol{u}(t)\rVert + \lVert\boldsymbol{v}(t; \boldsymbol{u})\rVert \leq 2 = v_{\max}.
    \end{multline*}
    The fixed-point iteration algorithm~\eqref{eq:seq_tau} solves the Problems~\ref{prob:1}--\ref{prob:2}. Note, that this algorithm defines a constructive way of finding a minimal root for Eq.~25 of~\cite{Bakolas2014-ze} in the problem of interception in the position space.

    The fixed-point iteration algorithm can be sped up by using the best universal lower estimator $T_{\rho}$ instead $\tau_{\rho}$ (see \cite{Buzikov2022-rb}). The best universal lower estimator provides the largest step $t_i - t_{i - 1}$ in the fixed-point iteration algorithm for the case when only the Lipschitz constant $v$ is known regarding the target trajectory $\boldsymbol{h}_T \in \mathrm{Lip}_v(\mathbb{R}^+_0, \mathcal{H})$ at the time of designing the algorithm. The best lower estimator is given by
    \begin{equation}\label{eq:best_est}
        T_\rho(t, \mathbf{h}) = \min\left\{\theta \in \mathbb{R}^+_0:\theta \geq t,\: \rho(\theta, \mathbf{h}) = v(\theta - t) + \ell\right\}.
    \end{equation}
    Here, as before, $\rho_{\mathcal{R}_{\boldsymbol{r}}}$ or $\rho_{\mathcal{R}_{\boldsymbol{v}}}$ is substituted in place of the function $\rho$, depending on the problem to be solved. The calculation of $T_\rho(t, \mathbf{h})$ requires knowledge of the solution of $\rho(\tau, \mathbf{h}) = v(\tau - t) + \ell$ which is hard to obtain in the general case. Further, we obtain the solution for an important special case of zero initial velocity of the isotropic rocket.

    \begin{thm}\label{th:Tr_calc}
        If $v_0 = 0$ and $\rho_{\mathcal{R}_{\boldsymbol{r}}}(t, \mathbf{h}) > \ell$, then the best universal lower estimator $T_{\rho_{\mathcal{R}_{\boldsymbol{r}}}}$ can be computed using the principle branch of Lambert $W_0$-function as follows:
        \begin{multline*}
            T_{\rho_{\mathcal{R}_{\boldsymbol{r}}}}(t, \mathbf{h}) = t + \frac{\lVert \mathbf{h} \rVert + 1 - \ell - t}{1 + v}\\
            + W_0\left(-\frac{e^{-\frac{1 + vt - \ell + \lVert \mathbf{h} \rVert}{1 + v}}}{1 + v}\right).
        \end{multline*}
    \end{thm}
    \begin{pf}
        According to~\eqref{eq:best_est}, the value of $T_{\rho_{\mathcal{R}_{\boldsymbol{r}}}}(t, \mathbf{h})$ can be calculated as the smallest root of the equation
        \begin{equation*}
            \lVert\mathbf{h}\rVert - (\theta - 1 + e^{-\theta}) = v(\theta - t) + \ell
        \end{equation*}
        with respect to $\theta \geq t$. Let's convert this equation to the form
        \begin{multline*}
            \left(\theta - \frac{1 + vt - \ell + \lVert \mathbf{h} \rVert}{1 + v}\right)e^{\theta - \frac{1 + vt - \ell + \lVert \mathbf{h} \rVert}{1 + v}}\\
            = -\frac{e^{-\frac{1 + vt - \ell + \lVert \mathbf{h} \rVert}{1 + v}}}{1 + v}.
        \end{multline*}
        Using the property of Lambert functions $W(x)e^{W(x)} = x$, we obtain
        \begin{equation*}
            \theta = t + \frac{\lVert \mathbf{h} \rVert + 1 - \ell - t}{1 + v} + W\left(-\frac{e^{-\frac{1 + vt - \ell + \lVert \mathbf{h} \rVert}{1 + v}}}{1 + v}\right).
        \end{equation*}
        Now it remains to show that we should use the $W_0$-branch of the multivalued Lambert $W$-function. Since for negative real arguments the Lambert function has only at most two real branches $W_0$ and $W_{-1}$(see \cite{Corless1996-go}), we must show that using the $W_{-1}$-branch is incorrect. According to Lemma 3.4 of \cite{Buzikov2022-rb},
        $T_{\rho_{\mathcal{R}_{\boldsymbol{r}}}}(t, \mathbf{h}) \geq \tau_{\rho_{\mathcal{R}_{\boldsymbol{r}}}}(t, \mathbf{h})$. Thus,
        \begin{multline}\label{eq:W_neq}
            t + \frac{\lVert\mathbf{h}\rVert + 1 - \ell - t}{1 + v} + W_{-1}\left(-\frac{e^{-\frac{1 + vt - \ell + \lVert\mathbf{h}\rVert}{1 + v}}}{1 + v}\right) \\
            - \left(t + \frac{\lVert\mathbf{h}\rVert - (t - 1 + e^{-t}) - \ell}{1 + v}\right)\\
            = W_{-1}\left(-\frac{e^{-\frac{1 + vt - \ell + \lVert\mathbf{h}\rVert}{1 + v}}}{1 + v}\right) + \frac{e^{-t}}{1 + v} \geq 0.
        \end{multline}
        The lower branch of the Lambert function satisfies the inequality $W_{-1}(x) \leq -1$ (see \cite{Corless1996-go}). Therefore,
       \begin{equation*}
            W_{-1}\left(-\frac{e^{-\frac{1 + vt - \ell + \lVert\mathbf{h}\rVert}{1 + v}}}{1 + v}\right) + \frac{e^{-t}}{1 + v} \leq -1 + \frac{e^{-t}}{1 + v} \leq 0.
        \end{equation*} 
        The last inequality converts to an equality only if $t = 0$ and $v = 0$. Hence, for $t > 0$ or $v > 0$, this inequality contradicts~\eqref{eq:W_neq}. Let us consider separately the case $t = 0$ and $v = 0$. Using~\eqref{eq:W_neq}, we obtain
        \begin{equation*}
            W_{-1}\left(-e^{-1 + \ell - \lVert\mathbf{h}\rVert}\right) = -1.
        \end{equation*}
        Using the fact that the function $W_{-1}$ reaches the value $-1$ only at the point $-1/e$, after some simple transformations we get $\ell = \lVert\mathbf{h}\rVert$, which contradicts $\rho_{\mathcal{R}_{\boldsymbol{r}}}(t, \mathbf{h}) > \ell$.\qed
    \end{pf}

    \begin{thm}\label{th:Tv_calc}
        If $v_0 = 0$ and $\rho_{\mathcal{R}_{\boldsymbol{v}}}(t, \mathbf{h}) > \ell$, then the best universal lower estimator $T_{\rho_{\mathcal{R}_{\boldsymbol{v}}}}$ can be computed using the principle branch of Lambert $W_0$-function as follows:
        \begin{equation*}
            T_{\rho_{\mathcal{R}_{\boldsymbol{v}}}}(t, \mathbf{h}) =
            \begin{cases}
                +\infty, \quad v = 0, \: 1 + \ell \leq \lVert\mathbf{h}\rVert;\\
                -\ln(1 + \ell - \lVert\mathbf{h}\rVert), \quad v = 0, \: 1 + \ell > \lVert\mathbf{h}\rVert;\\
                \frac{vt + \lVert\mathbf{h}\rVert - 1 - \ell}{v} + W_0\left(\frac{e^{-t + \frac{1 + \ell - \lVert\mathbf{h}\rVert}{v}}}{v}\right), v > 0.
            \end{cases}
        \end{equation*}
    \end{thm}
    \begin{pf}
        According to~\eqref{eq:best_est}, the value of $T_{\rho_{\mathcal{R}_{\boldsymbol{v}}}}(t, \mathbf{h})$ can be calculated as the smallest root of the equation $\lVert\mathbf{h}\rVert - (1 - e^{-\theta}) = v(\theta - t) + \ell$ with respect to $\theta \geq t$. We first consider separately the case of a resting target $v = 0$. It is not difficult to obtain a formal solution to this equation: $\theta = -\ln(1 + \ell - \lVert\mathbf{h}\rVert)$. If the argument of the logarithm is not positive, i.e., $1 + \ell \leq \lVert\mathbf{h}\rVert$, then this means that the original equation has no solution, so formally reaching the velocity $\mathbf{h}$ is impossible. This is to be expected, since for a given problem $\mathbf{h}$ means the required velocity, and an isotropic rocket cannot reach the speed $\lVert\mathbf{h}\rVert - \ell \geq 1$ in any way.
        
        Consider the case $v > 0$. Let's transform the equation to the form
        \begin{multline*}
            \left(\theta - t + \frac{1 + \ell - \lVert\mathbf{h}\rVert}{v}\right)e^{\theta - t + \frac{1 + \ell - \lVert\mathbf{h}\rVert}{v}}\\
            = \frac{1}{v}e^{- t + \frac{1 + \ell - \lVert\mathbf{h}\rVert}{v}} > 0.
        \end{multline*}
        Using the property of Lambert functions $W(x)e^{W(x)} = x$ and considering that for non-negative arguments there exists only the principal branch of the Lambert function $W_0$ (of the real-valued ones), we obtain
        \begin{equation*}
            \theta = t + \frac{\lVert\mathbf{h}\rVert - 1 - \ell}{v} + W_0\left(\frac{1}{v}e^{-t + \frac{1 + \ell - \lVert\mathbf{h}\rVert}{v}}\right).
        \end{equation*}\qed
    \end{pf}
    The above theorems give expressions for the best universal lower estimators for the important special case of zero initial velocity of the isotropic rocket. Using these functions, according to Theorem~3.3 of \cite{Buzikov2022-rb}, we can construct a sequence that converges to the minimal interception time $T^*[\boldsymbol{h}_T]$:
    \begin{equation}\label{eq:seq_T}
        t_i = T_\rho(t_{i - 1}, \boldsymbol{h}_T(t_{i - 1})), \quad t_0 = 0.
    \end{equation}
    Moreover, compared to the sequence~\eqref{eq:seq_tau}, the step $t_i - t_{i - 1}$ will always be no smaller than the same step of the sequence~\eqref{eq:seq_tau}, which speeds up convergence. Also, the sequence~\eqref{eq:seq_T} has the optimality property in the sense that among all monotone sequences based on universal lower estimators, the step of this sequence is the largest (see~\cite{Buzikov2022-rb}).

    In practice, sequences~\eqref{eq:seq_tau},~\eqref{eq:seq_T} can be used if the relative error in approaching the final position is fixed so that the number of iterations determined by sequences~\eqref{eq:seq_tau},~\eqref{eq:seq_T} is finite. Let $\varepsilon \in \mathbb{R}^+$ is the relative error, then we must compute $t_i$ using~\eqref{eq:seq_tau} or~\eqref{eq:seq_T} until $\rho(t_i, \boldsymbol{h}_T(t_i)) < \ell(1 + \varepsilon)$. If $k \in \mathbb{N}$ is the minimal number such that $\rho(t_k, \boldsymbol{h}_T(t_k)) < \ell(1 + \varepsilon)$, then we guarantee that the isotropic rocket can be transferred in time $t_k$ to the ball with the radius $\ell(1 + \varepsilon)$ and the center at $\boldsymbol{h}_T(t_k)$.

    \section{Simulation examples}

\begin{figure*}
        \centering
        \includegraphics[width=0.90\textwidth]{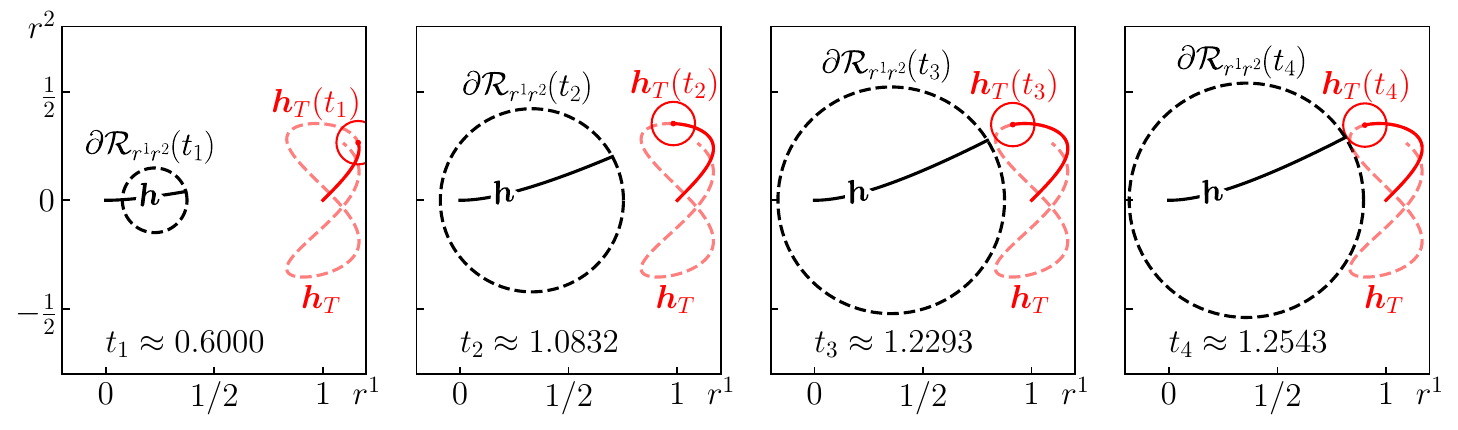}
        \caption{First steps of the fixed-point iteration algorithm $t_i = \tau_{\rho_{\mathcal{R}_{\boldsymbol{r}}}}(t_{i - 1}, \boldsymbol{h}_T(t_{i - 1}))$, $t_0 = 0$ for initial speed $v_0 = 1/2$ and $v_{\max} = 1$. The red circle is the boundary of the capture region with $\ell = 1/10$. The isotropic rocket must approach to the position $\boldsymbol{h}_T(t)$ (see~\eqref{eq:ex1}) at distance $\ell$ as fast as possible}
        \label{fig:r1r2tau}
    \end{figure*}

    In this section we illustrate the steps of the fixed-point iteration algorithm based on the universal lower estimators given by~\eqref{eq:seq_tau} and~\eqref{eq:seq_T}. As an example of Problem~\ref{prob:1}, we consider the motion of a target on a plane in the two-dimensional position space. Let the trajectory of the target be described by the following parameterization of the Lissajous curve:
    \begin{equation}\label{eq:ex1}
        \boldsymbol{h}_T(t) =
        \begin{bmatrix}
            1 + \frac16\sin3t\\
            \frac{\sqrt2}{4}\sin\sqrt2t
        \end{bmatrix}.
    \end{equation}
    This trajectory belongs to the class $\mathrm{Lip}_{1/2}(\mathbb{R}^+_0, \mathbb{R}^2)$, so $v = 1/2$. For non-zero initial velocity of the isotropic rocket, the first steps of the fixed-point iteration algorithm are shown in Fig.~\ref{fig:r1r2tau}.
    
    \begin{figure*}
        \centering
        \includegraphics[width=0.90\textwidth]{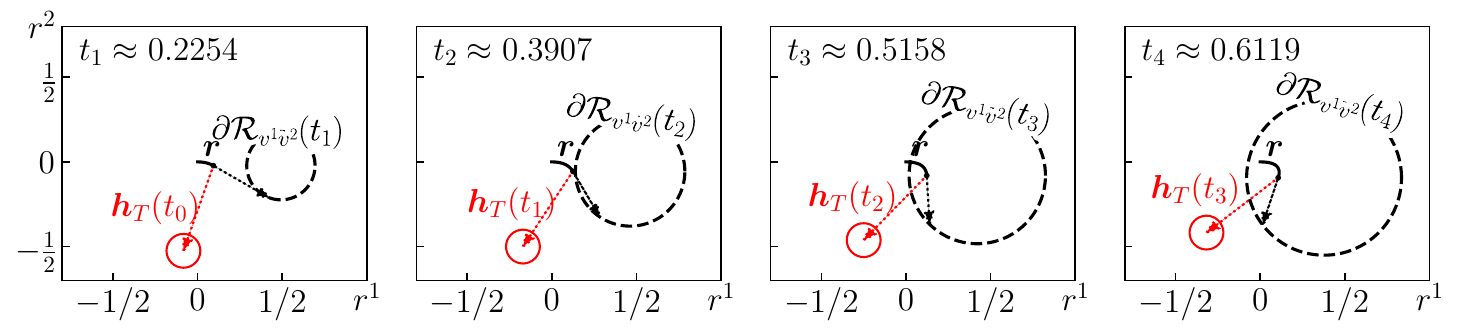}
        \caption{First steps of the fixed-point iteration algorithm $t_i = \tau_{\rho_{\mathcal{R}_{\boldsymbol{v}}}}(t_{i - 1}, \boldsymbol{h}_T(t_{i - 1}))$, $t_0 = 0$ for initial speed $v_0 = 1/2$ and $v_{\max} = 2$. The red circle is the boundary of the capture region with $\ell = 1/10$. The isotropic rocket must reach the velocity $\boldsymbol{h}_T(t)$ (see~\eqref{eq:ex2}) with precision $\ell$ as fast as possible. The black dotted arrow is a velocity of the isotropic rocket}
        \label{fig:v1v2tau}
    \end{figure*}

    For the Problem~\ref{prob:2}, we set $n = 2$ and
    \begin{equation}\label{eq:ex2}
        \boldsymbol{h}_T(t) =
        \begin{bmatrix}
            -\frac{8}{15}\sin\frac{3t}2\\
            -\frac{8}{15}\cos\frac{3t}2
        \end{bmatrix}.
    \end{equation}
    Here, $\boldsymbol{h}_T \in \mathrm{Lip}_{4/5}(\mathbb{R}^+_0, \mathbb{R}^2)$ defines the desired value of velocity changing on time (see Fig.~\ref{fig:v1v2tau}).
    
\begin{figure*}
        \centering
        \includegraphics[width=0.84\textwidth]{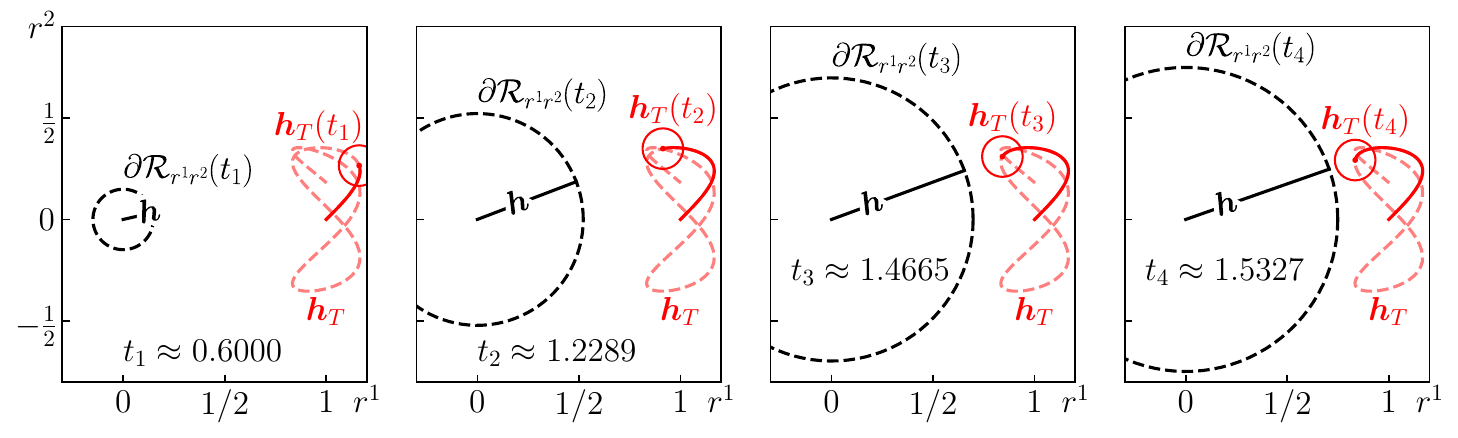}
        \includegraphics[width=0.84\textwidth]{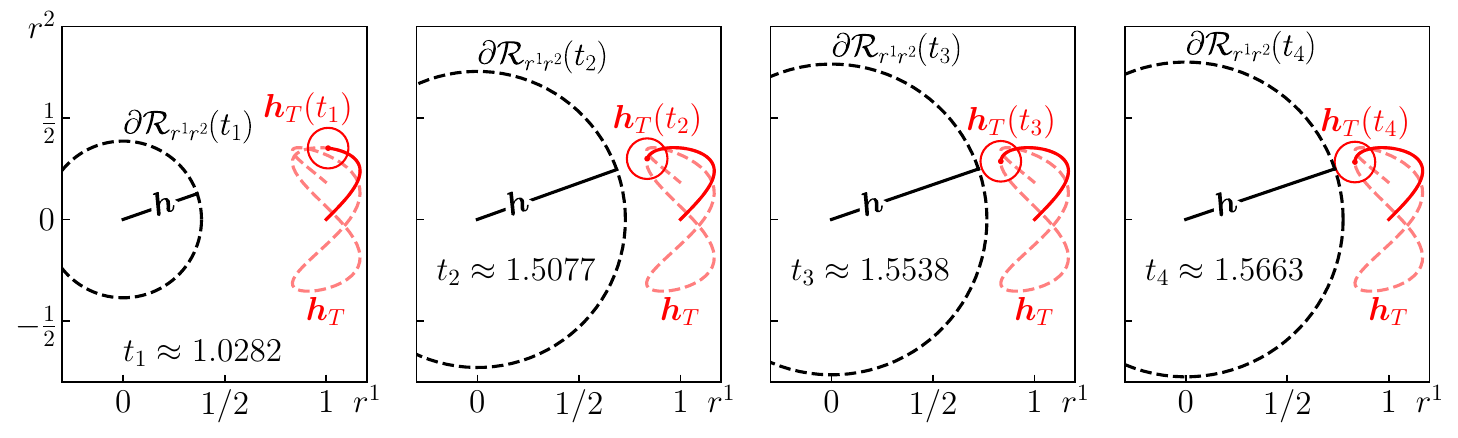}
        \caption{First steps of the fixed-point iteration algorithm based on $\tau_{\rho_{\mathcal{R}_{\boldsymbol{r}}}}$ (top row) and $T_{\rho_{\mathcal{R}_{\boldsymbol{r}}}}$ (bottom row). The initial speed is $v_0 = 0$ and $v_{\max} = 1$. The red circle is the boundary of the capture region with $\ell = 1/10$. The isotropic rocket must approach to the position $\boldsymbol{h}_T(t)$ (see~\eqref{eq:ex1}) at distance $\ell$ as fast as possible}
        \label{fig:tau_vs_T_r}
    \end{figure*}

    Due to the Theorem~\ref{th:Tr_calc}, it is possible to efficiently compute the best universal lower estimator for zero initial speed, so it is possible to compare the iterative progress with the algorithm based on the simple universal lower estimator. The first steps of these two methods are depicted in Fig~\ref{fig:tau_vs_T_r}. As expected, the method based on the function $T_{{\rho_{\mathcal{R}_{\boldsymbol{r}}}}}$ converges to the solution no slower than the method based on the function $\tau_{\rho_{\mathcal{R}_{\boldsymbol{r}}}}$.

    Taking into account that the rate of convergence of sequences~\eqref{eq:seq_tau},~\eqref{eq:seq_T} can be arbitrary slow (see~\cite{Buzikov2022-rb}), it is reasonable to use the described algorithm to obtain a good initial approximation for Newton’s method for finding the root of the equation $\rho(t, \boldsymbol{h}_T(t)) = \ell$.

    \section{Conclusions}
    
    In this study, an exhaustive description of the reachable set of an isotropic rocket is provided for the case of a non-zero drag coefficient. For this purpose, a parametric class of control inputs sufficient to reach any point on the boundary of the reachable set was identified using the maximum principle (Lemma~\ref{lem:ext_cont}). The expressions~\eqref{eq:vE} and~\eqref{eq:rE} describe extremal trajectories that terminate at the boundary of the reachable set. These expressions provide an explicit analytical description of the boundary of the reachable set~\eqref{eq:dR} for a given moment in time. Using this description, always convergent fixed-point iteration algorithms given by~\eqref{eq:seq_tau},~\eqref{eq:seq_T} are constructed for the problem of the fastest reaching of a moving position and for the problem of reaching a given time-varying velocity using an explicit description of the distance function from the corresponding projection of the reachable set to a given point. Theorems~\ref{th:Tr_calc}--\ref{th:Tv_calc} allow us to achieve the largest guaranteed step of the proposed algorithm for arbitrary Lipschitz motions of a target.

    A further development of this study is the construction of a similar numerical procedure for the problem of the fastest transfer to a given moving state (both the coordinate and velocity are important for intercepting). Computing the distance is a rather difficult problem because of the complex geometry of the reachable set. However, the properties of convexity and closeness make it possible to construct numerical procedures to estimate the distance to the reachable set.

    \begin{ack}
        The research was supported by RSF (project No. 23-19-00134).
    \end{ack}

    \appendix
    
    \section{Deriving extremal trajectories}\label{ap:ex_traj}
    
    We first consider the case $\lVert\boldsymbol{\lambda}_T\rVert\lVert\boldsymbol{\eta}_T\rVert > |(\boldsymbol{\lambda}_T, \boldsymbol{\eta}_T)|$, when the conditions of Lemma~\ref{lem:eta_zero} are not satisfied. Let $\boldsymbol{\xi}_T = \boldsymbol{\eta}_T - \boldsymbol{\lambda}_T \neq \boldsymbol{0}$. Then
    \begin{multline*}
        \int\limits_0^t \boldsymbol{u}_E(s; T, \mathbf{p}_T)e^{s - t}\mathrm{d}s
        = \int\limits_0^t \frac{\boldsymbol{\lambda}_T + \boldsymbol{\xi}_Te^{s - T}}{\lVert\boldsymbol{\lambda}_T + \boldsymbol{\xi}_Te^{s - T}\rVert} e^{s - t}\mathrm{d}s
    \end{multline*}
    \begin{multline*}
        = e^{T - t}\int\limits_{e^{-T}}^{e^{t - T}} \frac{\boldsymbol{\lambda}_T + \boldsymbol{\xi}_T\xi}{\lVert\boldsymbol{\lambda}_T + \boldsymbol{\xi}_T\xi\rVert}\mathrm{d}\xi\\
        = \boldsymbol{\lambda}_Te^{T - t} \int\limits_{e^{-T}}^{e^{t - T}} \frac{\mathrm{d}\xi}{\sqrt{\lVert\boldsymbol{\xi}_T\rVert^2\xi^2 + 2(\boldsymbol{\lambda}_T, \boldsymbol{\xi}_T)\xi + \lVert\boldsymbol{\lambda}_T\rVert^2}}\\
        + \boldsymbol{\xi}_Te^{T - t} \int\limits_{e^{-T}}^{e^{t - T}} \frac{\xi\mathrm{d}\xi}{\sqrt{\lVert\boldsymbol{\xi}_T\rVert^2\xi^2 + 2(\boldsymbol{\lambda}_T, \boldsymbol{\xi}_T)\xi + \lVert\boldsymbol{\lambda}_T\rVert^2}}.
    \end{multline*}
    Using the reference integrals 380.001 and 380.011 (see \cite{Dwight1957-tb}), we obtain
    \begin{multline*}
        \int\limits_{e^{-T}}^{e^{t - T}} \frac{\mathrm{d}\xi}{\sqrt{\lVert\boldsymbol{\xi}_T\rVert^2\xi^2 + 2(\boldsymbol{\lambda}_T, \boldsymbol{\xi}_T)\xi + \lVert\boldsymbol{\lambda}_T\rVert^2}}\\
        = \frac{1}{\lVert\boldsymbol{\xi}_T\rVert} \ln{\left(\lVert\boldsymbol{\xi}_T\rVert\lVert\boldsymbol{\lambda}_T + \boldsymbol{\xi}_T\xi\rVert + \rVert\boldsymbol{\xi}_T\lVert^2\xi + (\boldsymbol{\lambda}_T, \boldsymbol{\xi}_T)\right)}\Bigr|_{e^{-T}}^{e^{t - T}}\\
        = \frac{1}{\lVert\boldsymbol{\xi}_T\rVert} \ln{\frac{\lVert\boldsymbol{\xi}_T\rVert\lVert\boldsymbol{\eta}(t; T, \mathbf{p}_T)\rVert + (\boldsymbol{\xi}_T, \boldsymbol{\eta}(t; T, \mathbf{p}_T))}{\lVert\boldsymbol{\xi}_T\rVert\lVert\boldsymbol{\eta}(0; T, \mathbf{p}_T)\rVert + (\boldsymbol{\xi}_T, \boldsymbol{\eta}(0; T, \mathbf{p}_T))}};
    \end{multline*}
    \begin{multline*}
        \int\limits_{e^{-T}}^{e^{t - T}} \frac{\xi\mathrm{d}\xi}{\sqrt{\lVert\boldsymbol{\xi}_T\rVert^2\xi^2 + 2(\boldsymbol{\lambda}_T, \boldsymbol{\xi}_T)\xi + \lVert\boldsymbol{\lambda}_T\rVert^2}}\\
        = \frac{\lVert\boldsymbol{\lambda}_T + \boldsymbol{\xi}_T\xi\rVert}{\lVert\boldsymbol{\xi}_T\rVert^2}\Bigr|_{e^{-T}}^{e^{t - T}} \\
        - \frac{(\boldsymbol{\lambda}_T, \boldsymbol{\xi}_T)}{\lVert\boldsymbol{\xi}_T\rVert^2}\int\limits_{e^{-T}}^{e^{t - T}} \frac{\mathrm{d}\xi}{\sqrt{\lVert\boldsymbol{\xi}_T\rVert^2\xi^2 + 2(\boldsymbol{\lambda}_T, \boldsymbol{\xi}_T)\xi + \lVert\boldsymbol{\lambda}_T\rVert^2}}\\
        = \frac{\lVert\boldsymbol{\eta}(t; T, \mathbf{p}_T)\rVert - \lVert\boldsymbol{\eta}(0; T, \mathbf{p}_T)\rVert}{\lVert\boldsymbol{\xi}_T\rVert^2} \\
        -\frac{(\boldsymbol{\lambda}_T, \boldsymbol{\xi}_T)}{\lVert\boldsymbol{\xi}_T\rVert^3}\ln{\frac{\lVert\boldsymbol{\xi}_T\rVert\lVert\boldsymbol{\eta}(t; T, \mathbf{p}_T)\rVert + (\boldsymbol{\xi}_T, \boldsymbol{\eta}(t; T, \mathbf{p}_T))}{\lVert\boldsymbol{\xi}_T\rVert\lVert\boldsymbol{\eta}(0; T, \mathbf{p}_T)\rVert + (\boldsymbol{\xi}_T, \boldsymbol{\eta}(0; T, \mathbf{p}_T))}}.
    \end{multline*}
    Thus,
    \begin{multline}\label{eq:vE_1}
        \boldsymbol{v}_E(t; T, \mathbf{p}_T) = \frac{\lVert\boldsymbol{\eta}(t; T, \mathbf{p}_T)\rVert - \lVert\boldsymbol{\eta}(0; T, \mathbf{p}_T)\rVert}{\lVert\boldsymbol{\xi}_T\rVert^2}e^{T - t}\boldsymbol{\xi}_T\\
        + \boldsymbol{v}_0 e^{-t} + \left(\frac{\boldsymbol{\lambda}_T}{\lVert\boldsymbol{\xi}_T\rVert} - \frac{(\boldsymbol{\lambda}_T, \boldsymbol{\xi}_T)}{\lVert\boldsymbol{\xi}_T\rVert^3}\boldsymbol{\xi}_T\right)e^{T - t}\\
        \cdot\ln{\frac{\lVert\boldsymbol{\xi}_T\rVert\lVert\boldsymbol{\eta}(t; T, \mathbf{p}_T)\rVert + (\boldsymbol{\xi}_T, \boldsymbol{\eta}(t; T, \mathbf{p}_T))}{\lVert\boldsymbol{\xi}_T\rVert\lVert\boldsymbol{\eta}(0; T, \mathbf{p}_T)\rVert + (\boldsymbol{\xi}_T, \boldsymbol{\eta}(0; T, \mathbf{p}_T))}}.
    \end{multline}

    For $\lVert\boldsymbol{\lambda}_T\rVert\lVert\boldsymbol{\eta}_T\rVert > |(\boldsymbol{\lambda}_T, \boldsymbol{\eta}_T)|$ we also obtain
    \begin{equation*}
        \int\limits_0^t \boldsymbol{u}_E(s; T, \mathbf{p}_T)\mathrm{d}s
        = \int\limits_0^t \frac{\boldsymbol{\lambda}_T + \boldsymbol{\xi}_Te^{s - T}}{\lVert\boldsymbol{\lambda}_T + \boldsymbol{\xi}_Te^{s - T}\rVert}\mathrm{d}s
    \end{equation*}
    \begin{multline*}
        = \int\limits_{e^{-T}}^{e^{t - T}} \frac{\boldsymbol{\lambda}_T + \boldsymbol{\xi}_T\xi}{\lVert\boldsymbol{\lambda}_T + \boldsymbol{\xi}_T\xi\rVert}\frac{\mathrm{d}\xi}{\xi}\\
        = \boldsymbol{\lambda}_T \int\limits_{e^{-T}}^{e^{t - T}} \frac{\mathrm{d}\xi}{\xi\sqrt{\lVert\boldsymbol{\xi}_T\rVert^2\xi^2 + 2(\boldsymbol{\lambda}_T, \boldsymbol{\xi}_T)\xi + \lVert\boldsymbol{\lambda}_T\rVert^2}}\\
        + \boldsymbol{\xi}_T\int\limits_{e^{-T}}^{e^{t - T}} \frac{\mathrm{d}\xi}{\sqrt{\lVert\boldsymbol{\xi}_T\rVert^2\xi^2 + 2(\boldsymbol{\lambda}_T, \boldsymbol{\xi}_T)\xi + \lVert\boldsymbol{\lambda}_T\rVert^2}}.
    \end{multline*}
    Using the reference integral 380.111 (see \cite{Dwight1957-tb}), we obtain
    \begin{multline*}
        \int\limits_{e^{-T}}^{e^{t - T}} \frac{\mathrm{d}\xi}{\xi\sqrt{\lVert\boldsymbol{\xi}_T\rVert^2\xi^2 + 2(\boldsymbol{\lambda}_T, \boldsymbol{\xi}_T)\xi + \lVert\boldsymbol{\lambda}_T\rVert^2}}\\
        = \frac{1}{\lVert\boldsymbol{\lambda}_T\rVert}\ln{\frac{\xi}{\lVert\boldsymbol{\lambda}_T\rVert\lVert\boldsymbol{\lambda}_T + \boldsymbol{\xi}_T\xi\rVert + \xi(\boldsymbol{\lambda}_T, \boldsymbol{\xi}_T) + \lVert\boldsymbol{\lambda}_T\rVert^2}}\Bigr|_{e^{-T}}^{e^{t - T}}\\
        = \frac{1}{\lVert\boldsymbol{\lambda}_T\rVert}\left(1 + \ln\frac{\lVert\boldsymbol{\lambda}_T\rVert\lVert\boldsymbol{\eta}(0; T, \mathbf{p}_T)\rVert + (\boldsymbol{\lambda}_T, \boldsymbol{\eta}(0; T, \mathbf{p}_T))}{\lVert\boldsymbol{\lambda}_T\rVert\lVert\boldsymbol{\eta}(t; T, \mathbf{p}_T)\rVert + (\boldsymbol{\lambda}_T, \boldsymbol{\eta}(t; T, \mathbf{p}_T))}\right).
    \end{multline*}

    Thus,
    \begin{multline}\label{eq:rE_1}
        \boldsymbol{r}_E(t; T, \mathbf{p}_T) = \boldsymbol{v}_0 - \boldsymbol{v}_E(t; T, \mathbf{p}_T)\\
        + \frac{\boldsymbol{\lambda}_T}{\lVert\boldsymbol{\lambda}_T\rVert}\left(t + \ln\frac{\lVert\boldsymbol{\lambda}_T\rVert\lVert\boldsymbol{\eta}(0; T, \mathbf{p}_T)\rVert + (\boldsymbol{\lambda}_T, \boldsymbol{\eta}(0; T, \mathbf{p}_T))}{\lVert\boldsymbol{\lambda}_T\rVert\lVert\boldsymbol{\eta}(t; T, \mathbf{p}_T)\rVert + (\boldsymbol{\lambda}_T, \boldsymbol{\eta}(t; T, \mathbf{p}_T))}\right)\\
        + \frac{\boldsymbol{\xi}_T}{\lVert\boldsymbol{\xi}_T\rVert}\ln{\frac{\lVert\boldsymbol{\xi}_T\rVert\lVert\boldsymbol{\eta}(t; T, \mathbf{p}_T)\rVert + (\boldsymbol{\xi}_T, \boldsymbol{\eta}(t; T, \mathbf{p}_T))}{\lVert\boldsymbol{\xi}_T\rVert\lVert\boldsymbol{\eta}(0; T, \mathbf{p}_T)\rVert + (\boldsymbol{\xi}_T, \boldsymbol{\eta}(0; T, \mathbf{p}_T))}}.
    \end{multline}

    Now consider the case $\lVert\boldsymbol{\lambda}_T\rVert\lVert\boldsymbol{\eta}_T\rVert = -(\boldsymbol{\lambda}_T, \boldsymbol{\eta}_T)$ and $\boldsymbol{\lambda}_T \neq \boldsymbol{0}$. According to the Lemma~\ref{lem:eta_zero}, the integrand function may have singularities, so splitting the integral into a sum of integrals is no longer allowed. The considered conditions allow us to express $\boldsymbol{\eta}_T$ by $\boldsymbol{\lambda}_T$: $\boldsymbol{\eta}_T = -\boldsymbol{\lambda}_T\lVert\boldsymbol{\eta}_T\rVert/\lVert\boldsymbol{\lambda}_T\rVert$.
    \begin{multline*}
        \int\limits_0^t \boldsymbol{u}_E(s; T, \mathbf{p}_T)e^{s - t}\mathrm{d}s\\
        = \int\limits_0^t \frac{\boldsymbol{\lambda}_T + (\boldsymbol{\eta}_T - \boldsymbol{\lambda}_T)e^{s - T}}{\lVert\boldsymbol{\lambda}_T + (\boldsymbol{\eta}_T - \boldsymbol{\lambda}_T)e^{s - T}\rVert} e^{s - t}\mathrm{d}s\\
        = e^{T - t}\int\limits_{e^{-T}}^{e^{t - T}} \frac{\boldsymbol{\lambda}_T + (\boldsymbol{\eta}_T - \boldsymbol{\lambda}_T)\xi}{\lVert\boldsymbol{\lambda}_T + (\boldsymbol{\eta}_T - \boldsymbol{\lambda}_T)\xi\rVert}\mathrm{d}\xi\\
        = e^{T - t}\frac{\boldsymbol{\lambda}_T}{\lVert\boldsymbol{\lambda}_T\rVert}\int\limits_{e^{-T}}^{e^{t - T}}\frac{\lVert\boldsymbol{\lambda}_T\rVert - (\lVert\boldsymbol{\lambda}_T\rVert + \lVert\boldsymbol{\eta}_T\rVert)\xi}{|\lVert\boldsymbol{\lambda}_T\rVert - (\lVert\boldsymbol{\lambda}_T\rVert + \lVert\boldsymbol{\eta}_T\rVert)\xi|}\mathrm{d}\xi
    \end{multline*}
    \begin{multline*}
        = -e^{T - t}\frac{\boldsymbol{\lambda}_T}{\lVert\boldsymbol{\lambda}_T\rVert}\int\limits_{e^{-T}}^{e^{t - T}}\mathrm{sgn}\left(\xi - \frac{\lVert\boldsymbol{\lambda}_T\rVert}{\lVert\boldsymbol{\lambda}_T\rVert + \lVert\boldsymbol{\eta}_T\rVert}\right)\mathrm{d}\xi \\
        = \frac{\boldsymbol{\lambda}_T}{\lVert\boldsymbol{\lambda}_T\rVert}\left(\left| e^{-t} - e^{\theta(T, \mathbf{p}_T) - t}\right| - \left|1 - e^{\theta(T, \mathbf{p}_T) - t}\right|\right).
        % = e^{T - t}\frac{\boldsymbol{\lambda}_T}{\lVert\boldsymbol{\lambda}_T\rVert}\left(\left|e^{-T} - \frac{\lVert\boldsymbol{\lambda}_T\rVert}{\lVert\boldsymbol{\lambda}_T\rVert + \lVert\boldsymbol{\eta}_T\rVert}\right| - \left|e^{t - T} - \frac{\lVert\boldsymbol{\lambda}_T\rVert}{\lVert\boldsymbol{\lambda}_T\rVert + \lVert\boldsymbol{\eta}_T\rVert}\right|\right).
    \end{multline*}
    \begin{multline*}
        \int\limits_0^t \boldsymbol{u}_E(s; T, \mathbf{p}_T)\mathrm{d}s
        = \int\limits_0^t \frac{\boldsymbol{\lambda}_T + (\boldsymbol{\eta}_T - \boldsymbol{\lambda}_T)e^{s - T}}{\lVert\boldsymbol{\lambda}_T + (\boldsymbol{\eta}_T - \boldsymbol{\lambda}_T)e^{s - T}\rVert}\mathrm{d}s\\
        = \frac{\boldsymbol{\lambda}_T}{\lVert\boldsymbol{\lambda}_T\rVert}\int\limits_0^t\frac{\lVert\boldsymbol{\lambda}_T\rVert - (\lVert\boldsymbol{\lambda}_T\rVert + \lVert\boldsymbol{\eta}_T\rVert)e^{s - T}}{|\lVert\boldsymbol{\lambda}_T\rVert - (\lVert\boldsymbol{\lambda}_T\rVert + \lVert\boldsymbol{\eta}_T\rVert)e^{s - T}|}\mathrm{d}s\\
        = -\frac{\boldsymbol{\lambda}_T}{\lVert\boldsymbol{\lambda}_T\rVert}\int\limits_0^t\mathrm{sgn}\left(s - \left(T + \ln\frac{\lVert\boldsymbol{\lambda}_T\rVert}{\lVert\boldsymbol{\lambda}_T\rVert + \lVert\boldsymbol{\eta}_T\rVert}\right)\right)\mathrm{d}s\\
        = \frac{\boldsymbol{\lambda}_T}{\lVert\boldsymbol{\lambda}_T\rVert}\left(\left|\theta(T, \mathbf{p}_T)\right| - \left|t - \theta(T, \mathbf{p}_T)\right|\right).
    \end{multline*}
    Therefore,
    \begin{multline}\label{eq:vE_2}
        \boldsymbol{v}_E(t; T, \mathbf{p}_T) = \boldsymbol{v}_0 e^{-t}\\
        + \frac{\boldsymbol{\lambda}_T}{\lVert\boldsymbol{\lambda}_T\rVert}\left(\left| e^{-t} - e^{\theta(T, \mathbf{p}_T) - t}\right| - \left|1 - e^{\theta(T, \mathbf{p}_T) - t}\right|\right),
    \end{multline}
    \begin{multline}\label{eq:rE_2}
        \boldsymbol{r}_E(t; T, \mathbf{p}_T) = \boldsymbol{v}_0 - \boldsymbol{v}_E(t; T, \mathbf{p}_T)\\
        + \frac{\boldsymbol{\lambda}_T}{\lVert\boldsymbol{\lambda}_T\rVert}\left(\left|\theta(T, \mathbf{p}_T)\right| - \left|t - \theta(T, \mathbf{p}_T)\right|\right).
    \end{multline}
    If $\lVert\boldsymbol{\lambda}_T\rVert\lVert\boldsymbol{\eta}_T\rVert = (\boldsymbol{\lambda}_T, \boldsymbol{\eta}_T)$ and $\boldsymbol{\lambda}_T \neq \boldsymbol{0}$, then we use $\boldsymbol{\eta}_T = \boldsymbol{\lambda}_T\lVert\boldsymbol{\eta}_T\rVert/\lVert\boldsymbol{\lambda}_T\rVert$. In this case, Lemma~\ref{lem:eta_zero} states that almost everywhere the extremal control input is equal to a constant vector:
    \begin{multline*}
        \int\limits_0^t \boldsymbol{u}_E(s; T, \mathbf{p}_T)e^{s - t}\mathrm{d}s
        \\
        = \int\limits_0^t \frac{\boldsymbol{\lambda}_T + (\boldsymbol{\eta}_T - \boldsymbol{\lambda}_T)e^{s - T}}{\lVert\boldsymbol{\lambda}_T + (\boldsymbol{\eta}_T - \boldsymbol{\lambda}_T)e^{s - T}\rVert} e^{s - t}\mathrm{d}s\\
        = \frac{\boldsymbol{\lambda}_T}{\lVert\boldsymbol{\lambda}_T\rVert} \int\limits_0^t \frac{1 - e^{s - T} + \frac{\lVert\boldsymbol{\eta}_T\rVert}{\lVert\boldsymbol{\lambda}_T\rVert}e^{s - T}}{\left|1 - e^{s - T} + \frac{\lVert\boldsymbol{\eta}_T\rVert}{\lVert\boldsymbol{\lambda}_T\rVert}e^{s - T}\right|}e^{s - t}\mathrm{d}s\\
        = \frac{\boldsymbol{\lambda}_T}{\lVert\boldsymbol{\lambda}_T\rVert} \int\limits_0^t e^{s - t}\mathrm{d}s = \frac{\boldsymbol{\lambda}_T}{\lVert\boldsymbol{\lambda}_T\rVert}(1 - e^{-t});
    \end{multline*}
    \begin{multline*}
        \int\limits_0^t \boldsymbol{u}_E(s; T, \mathbf{p}_T)\mathrm{d}s = \int\limits_0^t \frac{\boldsymbol{\lambda}_T + (\boldsymbol{\eta}_T - \boldsymbol{\lambda}_T)e^{s - T}}{\lVert\boldsymbol{\lambda}_T + (\boldsymbol{\eta}_T - \boldsymbol{\lambda}_T)e^{s - T}\rVert} \mathrm{d}s\\
        = \frac{\boldsymbol{\lambda}_T}{\lVert\boldsymbol{\lambda}_T\rVert} \int\limits_0^t \frac{1 - e^{s - T} + \frac{\lVert\boldsymbol{\eta}_T\rVert}{\lVert\boldsymbol{\lambda}_T\rVert}e^{s - T}}{\left|1 - e^{s - T} + \frac{\lVert\boldsymbol{\eta}_T\rVert}{\lVert\boldsymbol{\lambda}_T\rVert}e^{s - T}\right|}\mathrm{d}s\\
        = \frac{\boldsymbol{\lambda}_T}{\lVert\boldsymbol{\lambda}_T\rVert} \int\limits_0^t\mathrm{d}s = \frac{\boldsymbol{\lambda}_T}{\lVert\boldsymbol{\lambda}_T\rVert}t.
    \end{multline*}
    Thus,
   \begin{equation}\label{eq:vE_3}
        \boldsymbol{v}_E(t; T, \mathbf{p}_T) = \boldsymbol{v}_0 e^{-t} + \frac{\boldsymbol{\lambda}_T}{\lVert\boldsymbol{\lambda}_T\rVert}(1 - e^{-t}),
    \end{equation}
    \begin{equation}\label{eq:rE_3}
        \boldsymbol{r}_E(t; T, \mathbf{p}_T) = \boldsymbol{v}_0 - \boldsymbol{v}_E(t; T, \mathbf{p}_T) + \frac{\boldsymbol{\lambda}_T}{\lVert\boldsymbol{\lambda}_T\rVert}t.
    \end{equation} 
    
    The last case is $\boldsymbol{\lambda}_T = \boldsymbol{0}$. Since $\boldsymbol{\eta}_T \neq \boldsymbol{0}$, we obtain
    \begin{equation}\label{eq:vE_4}
        \boldsymbol{v}_E(t; T, \mathbf{p}_T) = \boldsymbol{v}_0 e^{-t} + \frac{\boldsymbol{\eta}_T}{\lVert\boldsymbol{\eta}_T\rVert}(1 - e^{-t}),
    \end{equation}
    \begin{equation}\label{eq:rE_4}
        \boldsymbol{r}_E(t; T, \mathbf{p}_T) = \boldsymbol{v}_0 - \boldsymbol{v}_E(t; T, \mathbf{p}_T) + \frac{\boldsymbol{\eta}_T}{\lVert\boldsymbol{\eta}_T\rVert}t.
    \end{equation}

    \bibliographystyle{plainnat}
    \bibliography{main}

\begin{thebibliography}{40}
\providecommand{\natexlab}[1]{#1}
\providecommand{\url}[1]{\texttt{#1}}
\expandafter\ifx\csname urlstyle\endcsname\relax
  \providecommand{\doi}[1]{doi: #1}\else
  \providecommand{\doi}{doi: \begingroup \urlstyle{rm}\Url}\fi

\bibitem[Abeysiriwardena and Das(2018)]{Abeysiriwardena2018-wv}
S~Abeysiriwardena and T~Das.
\newblock Optimal control of a self-propelled particle in a fluid flow field.
\newblock In \emph{Proc. Am. Control Conf.}, pages 4135--4140. IEEE, 2018.
\newblock \doi{10.23919/ACC.2018.8430993}.

\bibitem[Abeysiriwardena and Das(2019)]{Abeysiriwardena2019-yg}
S~Abeysiriwardena and T~K Das.
\newblock Energy-optimal guidance of autonomous underwater vehicles under fluid-particle interaction.
\newblock \emph{J. Guid. Control Dyn.}, 42\penalty0 (6):\penalty0 1393--1401, 2019.
\newblock \doi{10.2514/1.G003695}.

\bibitem[Akulenko(1996)]{Akulenko1996-oq}
L~D Akulenko.
\newblock Control synthesis in the problem of the time-optimal intersection of a sphere.
\newblock \emph{J. Appl. Math. Mech.}, 60\penalty0 (5):\penalty0 717--727, 1996.
\newblock \doi{10.1016/S0021-8928(96)00091-3}.

\bibitem[Akulenko(2008)]{Akulenko2008-so}
L~D Akulenko.
\newblock The time-optimal transfer of a perturbed dynamical object to a given position.
\newblock \emph{J. Appl. Math. Mech.}, 72\penalty0 (2):\penalty0 136--143, 2008.
\newblock \doi{10.1016/j.jappmathmech.2008.04.012}.

\bibitem[Akulenko(2011)]{Akulenko2011-kr}
L~D Akulenko.
\newblock Time-optimal steering of an object moving in a viscous medium to a desired phase state.
\newblock \emph{J. Appl. Math. Mech.}, 75\penalty0 (5):\penalty0 534--538, 2011.
\newblock \doi{10.1016/j.jappmathmech.2011.11.007}.

\bibitem[Akulenko and Koshelev(2003)]{Akulenko2003-xd}
L~D Akulenko and A~P Koshelev.
\newblock Time-optimal steering of a dynamic object to a given position under the equality of the initial and final velocities.
\newblock \emph{J. Comput. Syst. Sci. Int.}, 42\penalty0 (6):\penalty0 921--928, 2003.

\bibitem[Akulenko and Koshelev(2005)]{Akulenko2005-iz}
L~D Akulenko and A~P Koshelev.
\newblock Time optimal return of a dynamic object with the required velocity.
\newblock \emph{Dokl. Math.}, 72\penalty0 (1):\penalty0 653--657, 2005.

\bibitem[Akulenko and Koshelev(2007)]{Akulenko2007-aw}
L~D Akulenko and A~P Koshelev.
\newblock Time-optimal steering of a point mass to a specified position with the required velocity.
\newblock \emph{J. Appl. Math. Mech.}, 71\penalty0 (2):\penalty0 200--207, 2007.
\newblock \doi{10.1016/j.jappmathmech.2007.06.003}.

\bibitem[Akulenko and Shmatkov(2007)]{Akulenko2007-rq}
L~D Akulenko and A~M Shmatkov.
\newblock Time-optimal crossing of a sphere in a viscous medium.
\newblock \emph{J. Comput. Syst. Sci. Int.}, 46:\penalty0 19--26, 2007.
\newblock \doi{10.1134/S1064230707010042}.

\bibitem[Akulenko and Shmatkov(2018)]{Akulenko2018-ck}
L~D Akulenko and A~M Shmatkov.
\newblock Transfer of a dynamic object onto the surface of an ellipsoid.
\newblock \emph{J Comput Syst Sci Int}, 57\penalty0 (1):\penalty0 63--71, 2018.
\newblock \doi{10.1134/S1064230718010021}.

\bibitem[Bakolas(2014{\natexlab{a}})]{Bakolas2014-jc}
E~Bakolas.
\newblock Minimum time control for a newtonian particle in a spatiotemporal flow field.
\newblock In \emph{Proc. Am. Control Conf.}, pages 2342--2347. IEEE, 2014{\natexlab{a}}.
\newblock \doi{10.1109/ACC.2014.6858860}.

\bibitem[Bakolas(2014{\natexlab{b}})]{Bakolas2014-ze}
E~Bakolas.
\newblock Optimal guidance of the isotropic rocket in the presence of wind.
\newblock \emph{J. Optim. Theory Appl.}, 162\penalty0 (3):\penalty0 954--974, 2014{\natexlab{b}}.
\newblock \doi{10.1007/s10957-013-0504-4}.

\bibitem[Bakolas and Marchidan(2016)]{Bakolas2016-ry}
E~Bakolas and A~Marchidan.
\newblock Time-optimal control of a self-propelled particle in a spatiotemporal flow field.
\newblock \emph{Int. J. Control}, 89\penalty0 (3):\penalty0 623--634, 2016.
\newblock \doi{10.1080/00207179.2015.1088965}.

\bibitem[Bernhard(1970)]{Bernhard1970-po}
P~Bernhard.
\newblock Linear pursuit-evasion games and the isotropic rocket.
\newblock Technical Report SUDAAR No. 413, Stanford University, 1970.

\bibitem[Bernhard(1972)]{Bernhard1972-gd}
P~Bernhard.
\newblock Corner conditions for differential games.
\newblock \emph{IFAC Proc. Vol.}, 5\penalty0 (1, Part 4):\penalty0 367--372, 1972.
\newblock \doi{10.1016/S1474-6670(17)68356-3}.

\bibitem[Botkin et~al.(2011)Botkin, Hoffmann, and Turova]{Botkin2011-vp}
N~Botkin, K-H Hoffmann, and V~L Turova.
\newblock Stable numerical schemes for solving {Hamilton–Jacobi–Bellman–Isaacs} equations.
\newblock \emph{SIAM J. Sci. Comput.}, 33\penalty0 (2):\penalty0 992--1007, 2011.
\newblock \doi{10.1137/100801068}.

\bibitem[Botkin et~al.(2013)Botkin, Hoffmann, Mayer, and Turova]{Botkin2013-ox}
N~Botkin, K-H Hoffmann, N~Mayer, and V~Turova.
\newblock Computation of value functions in nonlinear differential games with state constraints.
\newblock In \emph{System Modeling and Optimization}, pages 235--244, Berlin, Heidelberg, 2013. Springer.
\newblock \doi{10.1007/978-3-642-36062-6\_24}.

\bibitem[Buzikov(2022)]{Buzikov2022-rb}
M~Buzikov.
\newblock Computing the minimum-time interception of a moving target.
\newblock \emph{arXiv:2210.03439}, 2022.

\bibitem[Corless et~al.(1996)Corless, Gonnet, Hare, Jeffrey, and Knuth]{Corless1996-go}
R~M Corless, G~H Gonnet, D~E~G Hare, D~J Jeffrey, and D~E Knuth.
\newblock On the lambert {W} function.
\newblock \emph{Adv. Comput. Math.}, 5\penalty0 (1):\penalty0 329--359, 1996.
\newblock \doi{10.1007/bf02124750}.

\bibitem[Dwight(1957)]{Dwight1957-tb}
H~B Dwight.
\newblock \emph{Tables of Integrals and Other Mathematical Data}.
\newblock The Macmillan Company, New York, USA, third edition, 1957.

\bibitem[Friedman(1971)]{Friedman1971-jn}
A~Friedman.
\newblock Computation of saddle points for differential games of pursuit and evasion.
\newblock \emph{Arch. Ration. Mech. Anal.}, 40\penalty0 (2):\penalty0 79--119, 1971.
\newblock \doi{10.1007/BF00250316}.

\bibitem[Gutman et~al.(1987)Gutman, Esh, and Gefen]{Gutman1987-wz}
S~Gutman, M~Esh, and M~Gefen.
\newblock Simple linear pursuit-evasion games.
\newblock \emph{Comput. Math. Appl.}, 13\penalty0 (1):\penalty0 83--95, 1987.
\newblock \doi{10.1016/0898-1221(87)90095-2}.

\bibitem[Isaacs(1955)]{Isaacs1955-bn}
R~Isaacs.
\newblock Differential games {IV}: Mainly examples.
\newblock Technical Report RM-1486, RAND Corporation, 1955.

\bibitem[Isaacs(1965)]{Isaacs1965-nd}
R~Isaacs.
\newblock \emph{Differential Games: A Mathematical Theory with Applications to Warfare and Pursuit, Control and Optimization}.
\newblock John Wiley and Sons, Inc., New York, 1965.

\bibitem[Kai et~al.(2019)Kai, Hamel, and Samson]{Kai2019-pi}
J-M Kai, T~Hamel, and C~Samson.
\newblock A unified approach to fixed-wing aircraft path following guidance and control.
\newblock \emph{Automatica}, 108:\penalty0 108491, 2019.
\newblock \doi{10.1016/j.automatica.2019.07.004}.

\bibitem[Kumkov and Patsko(2014)]{Kumkov2014-bn}
S~S Kumkov and V~S Patsko.
\newblock Phenomenon of narrow throats of level sets of value function in differential games.
\newblock \emph{Contrib. Game Theory Manag.}, 7:\penalty0 161--180, 2014.

\bibitem[Lee and Markus(1967)]{Lee1967-vf}
E~B Lee and L~Markus.
\newblock \emph{Foundations of Optimal Control Theory}.
\newblock Wiley, New York, reprint edition 1986 with corrections edition, 1967.

\bibitem[Lewin and Olsder(1989)]{Lewin1989-ly}
J~Lewin and G~J Olsder.
\newblock The isotropic rocket—a surveillance evasion game.
\newblock \emph{Comput Math Appl}, 18\penalty0 (1):\penalty0 15--34, 1989.
\newblock \doi{10.1016/0898-1221(89)90121-1}.

\bibitem[Manchester and Savkin(2002)]{Manchester2002-sl}
I~R Manchester and A~V Savkin.
\newblock Circular navigation guidance law for precision missile/target engagements.
\newblock In \emph{Proceedings of the 41st {IEEE} Conference on Decision and Control, 2002.}, volume~2, pages 1287--1292 vol.2. IEEE, 2002.
\newblock \doi{10.1109/CDC.2002.1184692}.

\bibitem[Melikian(1973)]{Melikian1973-ee}
A~A Melikian.
\newblock On minimal observations in a game of encounter.
\newblock \emph{J. Appl. Math. Mech.}, 37\penalty0 (3):\penalty0 407--414, 1973.
\newblock \doi{10.1016/0021-8928(73)90084-1}.

\bibitem[Pontryagin and Mishchenko(1969)]{Pontryagin1969-dt}
L~S Pontryagin and E~F Mishchenko.
\newblock The problem of the escape of one controlled object from another.
\newblock \emph{Dokl. Akad. Nauk}, 189\penalty0 (4):\penalty0 721--723, 1969.

\bibitem[Pontryagin et~al.(1962)Pontryagin, Boltyanskii, Gamkrelidze, and Mischenko]{Pontryagin1962-lt}
L~S Pontryagin, V~G Boltyanskii, R~V Gamkrelidze, and E~F Mischenko.
\newblock \emph{The Mathematical Theory of Optimal Processes}.
\newblock John Wiley \& Sons, Inc, New York, London, 1962.

\bibitem[Selvakumar and Bakolas(2015)]{Selvakumar2015-pb}
J~Selvakumar and E~Bakolas.
\newblock Optimal guidance of the isotropic rocket in a partially uncertain flow field.
\newblock In \emph{Proc. Eur. Control Conf.}, pages 3328--3333. IEEE, 2015.
\newblock \doi{10.1109/ECC.2015.7331048}.

\bibitem[Selvakumar and Bakolas(2018)]{Selvakumar2018-np}
J~Selvakumar and E~Bakolas.
\newblock Robust time-optimal guidance in a partially uncertain time-varying flow-field.
\newblock \emph{J. Optim. Theory Appl.}, 179\penalty0 (1):\penalty0 240--264, 2018.
\newblock \doi{10.1007/s10957-018-1326-1}.

\bibitem[Stieber and Fügenschuh(2022)]{Stieber2022-et}
A~Stieber and A~Fügenschuh.
\newblock Dealing with time in the multiple traveling salespersons problem with moving targets.
\newblock \emph{Cent Eur J Oper Res}, 30\penalty0 (3):\penalty0 991--1017, 2022.
\newblock \doi{10.1007/s10100-020-00712-7}.

\bibitem[Tankasala et~al.(2022)Tankasala, Pehlivanturk, Bakolas, and Pryor]{Tankasala2022-hn}
S~Tankasala, C~Pehlivanturk, E~Bakolas, and M~Pryor.
\newblock Generating smooth near time-optimal trajectories for steering drones.
\newblock In \emph{European Control Conference}, pages 1484--1490, 2022.
\newblock \doi{10.23919/ECC55457.2022.9838451}.

\bibitem[Vnuchkov(1998)]{Vnuchkov1998-jc}
D~N Vnuchkov.
\newblock Time-optimal steering of a dynamic system with linear dissipation to a given terminal state.
\newblock \emph{J. Comput. Syst. Sci. Int.}, 37\penalty0 (3):\penalty0 399--404, 1998.

\bibitem[Wong(1967)]{Wong1967-lu}
R~E Wong.
\newblock Some aerospace differential games.
\newblock \emph{J. Spacecr. Rockets}, 4\penalty0 (11):\penalty0 1460--1465, 1967.
\newblock \doi{10.2514/3.29114}.

\bibitem[Yao and Cao(2020)]{Yao2020-ye}
W~Yao and M~Cao.
\newblock Path following control in {3D} using a vector field.
\newblock \emph{Automatica}, 117:\penalty0 108957, 2020.
\newblock \doi{10.1016/j.automatica.2020.108957}.

\bibitem[Zheng(2020)]{Zheng2020-hx}
Z~Zheng.
\newblock Moving path following control for a surface vessel with error constraint.
\newblock \emph{Automatica}, 118:\penalty0 109040, 2020.
\newblock \doi{10.1016/j.automatica.2020.109040}.

\end{thebibliography}
    
\end{document}